\documentclass[10pt]{elsarticle}
\usepackage{latexsym,amscd,amssymb,epsfig, enumitem, amsmath, tikz, bm} 
\usepackage{booktabs}  
\usepackage{amsfonts}
\usepackage{amstext}
\usepackage{amsthm}
\usepackage{textcomp}  
\usepackage{xcolor}
\usepackage{amsmath}
\usepackage{tikz}
\usepackage{bm}
\usepackage{amsthm}
\usepackage{lscape} 
\usepackage{xargs} 
\usetikzlibrary{arrows.meta}
\usetikzlibrary{calc}
\usetikzlibrary{decorations.pathreplacing}
\usetikzlibrary{arrows.meta,calc}
\usepackage[colorinlistoftodos]{todonotes}

\theoremstyle{plain}
  \newtheorem{thm}{Theorem}[section]
  \newtheorem{prop}[thm]{Proposition}
  \newtheorem{lem}[thm]{Lemma}
  \newtheorem{cor}[thm]{Corollary}
  
\theoremstyle{definition}
  \newtheorem{defn}[thm]{Definition}
  \newtheorem{ex}[thm]{Example}

  \newdefinition{rmk}{Remark}
  \newproof{pf}{Proof}

\numberwithin{equation}{section}

\begin{document}

\title{The Edge-Product Space of Phylogenetic Trees is not Shellable}

\author{Grace Stadnyk}
\address{Department of Mathematics, Furman University, Greenville, SC 29613} 
\ead{grace.stadnyk@furman.edu}

\begin{abstract} The edge-product space of phylogenetic trees is a regular CW complex whose maximal closed cells correspond to trivalent trees with leaves labeled by a finite set $X$. The face poset of this cell decomposition is isomorphic to the Tuffley poset, a poset of labeled forests, with a unique minimum adjoined. 
We show that the edge-product space of phylogenetic trees is gallery-connected. We then use combinatorial properties of the Tuffley poset and a related graph known as NNI-tree space to show that, although open intervals of the Tuffley poset were proven to be shellable by Gill, Linusson, Moulton, and Steel, the edge-product space is not shellable. 

\end{abstract}


\maketitle 

\section{Introduction}

The edge-product space of phylogenetic trees, $\mathcal{E}(X)$, is a topological space arising from edge-weighted trees with leaves labeled by a finite set $X$. Broadly, $\mathcal{E}(X)$ is constructed via a map that associates each edge-weighted tree to a point in space whose coordinates are given by products of weights on edges between leaves.  This space is of interest in part because of its applications to molecular evolutionary biology and tree reconstruction using Markov processes (see \cite{gill}, \cite{semplesteel}), but our focus is on further studying its topological properties. A similar space of phylogenetic trees, known as BHV space, was studied earlier by Billera, Holmes, and Vogtmann in \cite{bhv}, the primary difference being that coordinates of points in BHV space correspond to sums of edge-weights, as opposed to products as in $\mathcal{E}(X)$. As noted in \cite{gill}, the edge-product space can be viewed as the compactification of the tree-space proposed in \cite{bhv}.

The combinatorics and topology of the edge-product space were first investigated by Moulton and Steel in \cite{oranges}, where $\mathcal{E}(X)$ was proven to have a CW complex decomposition with a face poset (omitting the empty face) that is isomorphic to a poset called the Tuffley poset. We will denote the Tuffley poset by $S(X)$. 
Gill, Linusson, Moulton, and Steel subsequently proved in \cite{gill} that the CW decomposition defined in \cite{oranges} is regular. This implies that the combinatorics of $S(X)$ encodes the topology of $\mathcal{E}(X)$. 

The regularity of the cell decomposition was proven in part by showing that intervals in the Tuffley poset (with a unique minimum adjoined) are shellable, which in turn implies that the order complex of such intervals are topological spheres. The shellability of intervals in $S(X) \cup \hat{0}$ was shown in \cite{gill} by establishing the existence of (but not an explicit construction of) recursive coatom orderings for these intervals. More recently in \cite{elect}, Hersh and Kenyon provided an explicit shelling of these intervals via an EC-labeling of the uncrossing poset, which has  wire diagrams that are in bijection with the elements of the Tuffley poset as its elements. 

We expand on the shellability results of Gill, Linusson, Moulton, and Steel with the following, which is the primary result of this paper: 

\begin{thm} \label{s5norco} For $|X| \geq 5$, the edge-product space $\mathcal{E}(X)$ is not shellable. 
\end{thm}
 
While the shellability of cell complexes was first explored in the context of simplicial and polyhedral complexes, we examine the shellability of the CW decomposition of $\mathcal{E}(X)$  in the sense of Bj{\"o}rner in \cite{bjorner}, who extended the notion to pure, regular CW complexes in a way that preserved many of the topological implications of shellability (see Definition \ref{dcwshell}). For example, a pure, shellable, regular CW complex has the homotopy type of a wedge of spheres and has shellable barycentric subdivision \cite{bjorner}. As is the case for polyhedral complexes, a convenient way to establish the shellability of a regular CW decomposition is to examine the combinatorics of its face poset. Specifically, Bj{\"o}rner proved that the face poset of a pure, regular CW complex augmented with a unique maximum is CL-shellable precisely when the CW complex is shellable (see Proposition \ref{shelliff}).

It is natural to ask whether the augmented face poset of $\mathcal{E}(X)$ is CL-shellable (thereby ascertaining whether the cell decomposition of $\mathcal{E}(X)$ is shellable) in part because several properties of the edge-product space and its face poset might lead one to suspect that it is. For one, the CL-shellability of intervals in the face poset of $\mathcal{E}(X)$ satisfies a prerequisite for CL-shellability of the face poset of $\mathcal{E}(X)$ in its entirety. In addition, we prove in Proposition \ref{galconnectedthm} that the edge-product space is gallery-connected, a property held by all shellable, regular CW complexes. It is also known that the augmented face poset of $\mathcal{E}(X)$ is CL-shellable when $|X| \leq 4$. 

As was the case in \cite{gill}, the primary method we use to establish Theorem \ref{s5norco} will be recursive coatom orderings, a powerful but technical tool first introduced by Bj{\"o}rner and Wachs in \cite{bw} to prove that a poset is lexicographically shellable and thus that its order complex is shellable. 
We show that a recursive coatom ordering of $\hat{S}(X)$, the face poset of $\mathcal{E}(X)$ augmented with a $\hat{1}$, cannot exist because of fundamental characteristics of a graph called the NNI-tree space.
As this graph is the facet-ridge graph of the cell decomposition of $\mathcal{E}(X)$, finding a recursive coatom ordering of $\hat{S}(X)$ is equivalent to finding a Hamiltonian path in NNI-tree space that satisfies certain additional properties. The existence of both cycles of length $l>3$ and connected, tree-like subgraphs within NNI-tree space obstruct, in some sense, the construction of such Hamiltonian paths and in turn the existence of recursive coatom orderings for $\hat{S}(X)$.

The non-shellability of $\mathcal{E}(X)$ fundamentally implies that no matter the order in which we attach the facets, or $d-$dimensional cells, of the CW complex as we build up $\mathcal{E}(X)$, we will always ultimately attach a facet whose intersectection with previously attached facets is not a pure, $(d-1)$-dimensional CW complex. Though $\mathcal{E}(X)$ is not shellable, this does not preclude the topological results that follow immediately from shellability. In particular, it is worth noting that Theorem \ref{s5norco} does not necessarily imply that the poset $S(X)$ is not shellable or, equivalently, that the barycentric subdivision of $\mathcal{E}(X)$ is not shellable. In fact, Vince and Wachs in \cite{vincewachs}, Walker in \cite{walker}, and Vince in \cite{vince} produced examples of regular CW complexes that are not shellable but that have shellable barycentric subdivisions. However, in this case, it indeed turns out that $S(X)$ is not shellable for $|X| \geq 5$; the intervals defined in the proof of Theorem \ref{s5norco} that prevent the existence of a recursive coatom ordering of $\hat{S}(X)$ also prevent the existence of a shelling of the order complex of $S(X)$.   

The edge-product space is an example of a toric cube, that is, the image of a cube under a map whose coordinates are given by monomials. In \cite{cubes},  Engstr{\"o}m, Hersh, and Sturmfels provide CW decompositions of these spaces. A potentially interesting avenue of future research, albeit one that moves away from phylogenetic applications, is to determine whether there exist other interesting and important classes of toric cubes that are shellable. 

After providing the necessary background on the edge-product space of phylogenetic trees, CW complex shellability, and NNI-tree space in Section \ref{background}, we will prove that the edge-product space is gallery-connected. We then turn to proving Theorem \ref{s5norco} in Section \ref{noshell}.

\section{Preliminaries}
\label{background}

\subsection{The Edge-product Space of Phylogenetic Trees} 
\label{epreview} 
We begin with a very brief overview of the edge-product space of phylogenetic trees and its face poset, known as the Tuffley poset. For more detailed information on this space, see \cite{oranges} and \cite{gill}.

Let $T$ be a tree with vertex set $V(T)$ and edge set $E(T)$. Let $\mathbf{T}(X)$ be the set of all trees with leaves labeled bijectively by the finite set $X$.  
 For any map $\lambda: E(T) \rightarrow [0,1]$, denote by $\lambda(E(T))$ the vector in $[0,1]^{E(T)}$ whose $e_i$ entry is $\lambda(e_i)$. We will call $\lambda(E(T))$ an \textit{edge-weight vector}. For $x, y \in X$ and a tree $T \in \mathbf{T}(X)$, let $p_{xy}$ denote the set of edges in the unique path between the leaf labeled $x$ and the leaf labeled $y$ in $T$.


Given a map $\lambda: E(T) \rightarrow [0,1]$ and a tree $T \in \mathbf{T}(X)$, define the map $f_{(T, \lambda)}$ as 
 \begin{align*}
f_{(T, \lambda)}: {X \choose 2} &\rightarrow [0, 1]  \nonumber \\
(x, y) &\mapsto \prod\limits_{e \in p_{xy}} \lambda(e)
\end{align*}
Let $f_{(T, \lambda)} {X \choose 2}$ be the vector in $[0,1]^{X \choose 2}$ whose $xy$ entry is $f_{(T, \lambda)}((x,y))$. For any edge-weight vector $\lambda$, call $f_{(T, \lambda)}{X \choose 2}$ an \textit{edge-product vector}. We then define $\Lambda_{T}$ as the map that sends each edge-weight vector to the corresponding edge-product vector: 
		\begin{align*} 
			\Lambda_{T}: [0, 1]^{E(T)} &\rightarrow [0, 1]^{X \choose 2} \nonumber 			\\
			 \lambda(E(T)) &\mapsto f_{(T, \lambda)} {X \choose 2} 
		\end{align*}

\begin{defn} Let $\mathcal{E}(X, T)$ be the image of the map $\Lambda_{T}$. The \textit{edge-product space for trees on $X$} is 
\begin{displaymath}
\mathcal{E}(X)=\bigcup_{T \in \mathbf{T}(X)} \mathcal{E}(X, T)
\end{displaymath}
\label{epspace}
\end{defn}

%

As noted above, one of the primary results of \cite{oranges} is that the edge-product space has a CW decomposition whose face poset is isomorphic to a poset called the Tuffley poset when we adjoin a unique minimal element. We will assume familiarity with the basics of CW complexes, though readers can refer to \cite{bjornertopmethods} and others for this as necessary. Before defining the Tuffley poset, however, let us review some necessary basic terminology of posets (partially ordered sets). 

We will refer to a poset $(P, \leq)$ as $P$ when no confusion could arise in doing so. For $x, y \in P$, we say $y$ \textit{covers} $x$ (equivalently $x$ is \textit{covered by} $y$) precisely when $x < y$ and no $z$ exists such that $x < z < y$. We denote this relationship, called a \textit{cover relation}, by $x \lessdot y$. A \textit{chain} in $P$ is a totally ordered subset of $P$. The \textit{length} of a chain is one less than the number of elements in the chain. A poset is \textit{graded} if all maximal chains have the same length. A poset is \textit{bounded} provided there exists a unique maximum in $P$, denoted $\hat{1}$, and a unique minimum in $P$, denoted $\hat{0}$. A \textit{coatom} of $P$ is any element of $P$ covered by $\hat{1}$. A poset is \textit{thin} if any interval containing a maximal chain of length two contains exactly 4 elements. Given an element $x \in P$, the \textit{principal order ideal generated by x}, denoted $I(x)$, is the induced subposet of $P$ consisting of all $y \leq x$ in $P$. Given a poset $P$, we can define an abstract simplicial complex called the \textit{order complex of $P$}, denoted $\Delta(P)$ that is defined by taking the $k$-faces to be the chains of length $k$ in $P$. When we say a poset $P$ is shellable, we mean that the order complex $\Delta(P)$ is shellable. See  \cite{stanley} or \cite{wachs} for further detail on posets.

Now we turn to defining the Tuffley poset, a primary focus of our study. The Tuffley poset has $X$-forests as elements and cover relations given by operations on $X$-forests called edge contraction and safe edge deletion. 

\begin{defn}
An \textit{$X$-tree} $\mathcal{T}$ is a pair $(T, \phi)$ consisting of a tree $T$ and a map $\phi: X \rightarrow V(T)$ such that all vertices in $V(T)-\phi(X)$ are of degree greater than 2. An \textit{$X$-forest} is a set $\mathcal{F}=\{(A, \mathcal{T}_A) : A \in \pi\}$ where $\pi$ is a set partition of $X$ and $\mathcal{T}_A=(T_A, \phi_A)$ is an $A$-tree for every block $A \in \pi$. For $v \in V(T)$, we say $v$ is \textit{labeled} if $v=\phi(x)$ for some $x \in X$. Otherwise $v$ is \textit{unlabeled}.  
\end{defn} 

\begin{defn} Let $e=(u, v)$ be an edge in the $X$-tree $\mathcal{T}_A = (T_A, \phi_A) \in \mathcal{F}$. A \textit{contraction} of the edge $e$ is the elimination of the edge $e$ in $\mathcal{F}$ with the identification of the vertices $u$ and $v$. The resulting vertex is labeled $\phi^{-1}_A(u) \cup \phi^{-1}_A(v)$, i.e. the resulting vertex is labeled by the union of labels on $u$ and $v$ in $\mathcal{F}$. A \textit{deletion} of the edge $e$ is the elimination of edge $e$ in $\mathcal{F}$ with no changes to the vertex set of $\mathcal{T}_A$ for any $\mathcal{T}_A \in \mathcal{F}$.  The deletion of edge $e=[u,v]$ is called a \textit{safe edge deletion} if both $u$ and $v$ are either labeled or of degree greater than 3.
 \label{deletion}
\end{defn}

\begin{defn} Let $\mathcal{F}=\{(A, \mathcal{T}_A): A \in \pi\}$ and $\mathcal{F'}=\{(B, \mathcal{T}_B): B \in \pi'\}$ be $X$-forests. Define a partial order $\leq$ on $X$-forests so that $\mathcal{F'} \leq \mathcal{F}$ if $\mathcal{F'}$ can be obtained from a sequence of contractions and safe deletions of edges of $\mathcal{F}$. The poset $(S(X), \leq)$ is called the \textit{Tuffley poset}.
\label{tuffleyposet}
\end{defn}

For a CW complex $\mathcal{K}$, there exists a related poset called the \textit{face poset} of $\mathcal{K}$ consisting of the set of closed cells of $\mathcal{K}$ ordered by containment with a unique minimum, $\hat{0}$, adjoined. We will denote the face poset of $\mathcal{K}$ by $P(\mathcal{K})$. Note that in \cite{gill} and \cite{oranges}, the face poset of the edge-product space $\mathcal{E}(X)$ is defined without the addition of a $\hat{0}$. We instead follow the convention of Bj{\"o}rner established in \cite{bjorner} by defining the face poset so that it contains a unique minimum, $\hat{0}$, representing the empty face.

Moulton and Steel proved that the Tuffley poset augmented with a unique minimum is isomorphic to the face poset of the edge-product space: 

\begin{thm} [{Theorem~3.3} \cite{oranges}] The edge-product space $\mathcal{E}(X)$ is a finite CW complex. 
The Tuffley poset with a unique minimum adjoined, $(S(X) \cup \hat{0}, \leq)$, is isomorphic to the face poset of $\mathcal{E}(X)$. 
\label{tuffleycw}
\end{thm}

The Tuffley poset has several nice properties; for example, it is thin and graded. The following are two properties of $S(X)$ that we will find especially useful for the results in this paper.

\begin{thm} [{Theorem~4.2} \cite{oranges}] Let $X$ be a finite set and $\mathcal{F}, \mathcal{F'} \in S(X)$.
\begin{enumerate} 
\item The element $\mathcal{F}$ is a maximal element of $S(X)$ if and only if $\mathcal{F}$ is a tree whose leaves are labeled bijectively by $X$ and whose internal vertices all have degree three.  
\item  The element $\mathcal{F'}$ is covered by $\mathcal{F}$ if and only if $\mathcal{F'}$ can be obtained from $\mathcal{F}$ by a single edge contraction or safe edge deletion. 
\end{enumerate}
\label{tuffleyprops}
\end{thm}

For the remainder of this paper we will denote the set of maximal elements of $S(X)$ by $\mathcal{C}_X$. 

\begin{figure}
\begin{center}
\begin{tikzpicture}[scale=0.6]
\node (1) at (0,0) {\begin{tikzpicture} [black, scale=0.2] 
	\node [below] at (0,0) {\footnotesize 1,2,3}; 
	\draw [fill] (0,0) circle [radius=0.1]; 	
	\end{tikzpicture}};
\node (2) at (4,0) {\begin{tikzpicture} [black, scale=0.2] 
	\node [below] at (0,0) {\footnotesize 1,2}; 
	\draw [fill] (0,0) circle [radius=0.1]; 
	\node [below] at (3,0) {\footnotesize 3};
	\draw [fill] (3,0) circle [radius=0.1];
	\end{tikzpicture}};
\node (3) at (8,0) {\begin{tikzpicture} [black, scale=0.2] 
	\node [below] at (0,0) {\footnotesize 1,3}; 
	\draw [fill] (0,0) circle [radius=0.1]; 
	\node [below] at (3,0) {\footnotesize 2};
	\draw [fill] (3,0) circle [radius=0.1];
	\end{tikzpicture}};
\node (4) at (12,0) {\begin{tikzpicture} [black, scale=0.2] 
	\node [below] at (0,0) {\footnotesize 1}; 
	\draw [fill] (0,0) circle [radius=0.1]; 
	\node [below] at (3,0) {\footnotesize 2,3};
	\draw [fill] (3,0) circle [radius=0.1];
	\end{tikzpicture}};
\node (5) at (16,0) {\begin{tikzpicture} [black, scale=0.2] 
	\node [below] at (0,0) {\footnotesize 1}; 
	\draw [fill] (0,0) circle [radius=0.1]; 
	\node [below] at (2,0) {\footnotesize 2};
	\draw [fill] (2,0) circle [radius=0.1];
	\node [below] at (4,0) {\footnotesize 3};
	\draw [fill] (4,0) circle [radius=0.1];
	\end{tikzpicture}};
\node (6) at (-2,4){\begin{tikzpicture} [black, scale=0.2]
	\node [below] at (0,0) {\footnotesize 3};
	\draw [fill] (0,0) circle [radius=0.1];
	\node [above] at (1,1) {\footnotesize 1,2};
	\draw [fill] (1,1) circle [radius=0.1];
	\draw (0,0)--(1,1);
	\end{tikzpicture}};
\node (7) at (2, 4){\begin{tikzpicture} [black, scale=0.2]
	\node [below] at (0,0) {\footnotesize 2};
	\draw [fill] (0,0) circle [radius=0.1];
	\node [above] at (1,1) {\footnotesize 1};
	\draw [fill] (1,1) circle [radius=0.1];
	\node [below] at (2,0) {\footnotesize 3};
	\draw [fill] (2,0) circle [radius=0.1];
	\draw (0,0)--(1,1);
	\end{tikzpicture}};
\node (8) at (6,4){\begin{tikzpicture} [black, scale=0.2]
	\node [above] at (0,1) {\footnotesize 1,3};
	\draw [fill] (0,1) circle [radius=0.1];
	\node [below] at (1,0) {\footnotesize 2};
	\draw [fill] (1,0) circle [radius=0.1];
	\draw (0,1)--(1,0);
	\end{tikzpicture}};
\node (9) at (10,4){\begin{tikzpicture} [black, scale=0.2]
	\node [below] at (0,0) {\footnotesize 2};
	\draw [fill] (0,0) circle [radius=0.1];
	\node [above] at (1,1) {\footnotesize 1};
	\draw [fill] (1,1) circle [radius=0.1];
	\node [below] at (2,0) {\footnotesize 3};
	\draw [fill] (2,0) circle [radius=0.1];
	\draw (2,0)--(1,1);
	\end{tikzpicture}};
\node (10) at (14,4){\begin{tikzpicture} [black, scale=0.2]
	\node [below] at (1,1) {\footnotesize 2,3};
	\draw [fill] (1,1) circle [radius=0.1];
	\node [above] at (1,2) {\footnotesize 1};
	\draw [fill] (1,2) circle [radius=0.1];
	\draw (1,1)--(1,2);
	\end{tikzpicture}};
\node (11) at (18,4){\begin{tikzpicture} [black, scale=0.2]
	\node [left] at (1,1) {\footnotesize 2};
	\draw [fill] (1,1) circle [radius=0.1];
	\node [below] at (2,0) {\footnotesize 3};
	\draw [fill] (2,0) circle [radius=0.1];
	\node [above] at (1,2) {\footnotesize 1};
	\draw [fill] (1,2) circle [radius=0.1];
	\draw (1,1)--(2,0);
	\end{tikzpicture}};
\node (12) at (0,8){\begin{tikzpicture} [black, scale=0.2]
	\node [below] at (0,0) {\footnotesize 2};
	\draw [fill] (0,0) circle [radius=0.1];
	\node [below] at (2,0) {\footnotesize 3};
	\draw [fill] (2,0) circle [radius=0.1];
	\node [above] at (1,1) {\footnotesize 1};
	\draw [fill] (1,1) circle [radius=0.1];
	\draw (0,0)--(1,1)--(2,0);
	\end{tikzpicture}};
\node (13) at (8,8){\begin{tikzpicture} [black, scale=0.2]
	\node [below left] at (1,1) {\footnotesize 2};
	\draw [fill] (1,1) circle [radius=0.1];
	\node [below] at (2,0) {\footnotesize 3};
	\draw [fill] (2,0) circle [radius=0.1];
	\node [above] at (1,2) {\footnotesize 1};
	\draw [fill] (1,2) circle [radius=0.1];
	\draw (1,1)--(1,2);
	\draw (1,1)--(2,0);
	\end{tikzpicture}};
\node (14) at (16,8){\begin{tikzpicture} [black, scale=0.2]
	\node [below] at (0,0) {\footnotesize 2};
	\draw [fill] (0,0) circle [radius=0.1];
	\node [below right] at (1,1) {\footnotesize 3};
	\draw [fill] (1,1) circle [radius=0.1];
	\node [above] at (1,2) {\footnotesize 1};
	\draw [fill] (1,2) circle [radius=0.1];
	\draw (0,0)--(1,1)--(1,2);
	\end{tikzpicture}};
\node (15) at (8,12){\begin{tikzpicture} [black, scale=0.2]
	\node [below] at (0,0) {\footnotesize 2};
	\draw [fill] (0,0) circle [radius=0.1];
	\node [above] at (1,1) {};
	\draw [fill] (1,1) circle [radius=0.1];
	\node [below] at (2,0) {\footnotesize 3};
	\draw [fill] (2,0) circle [radius=0.1];
	\node [above] at (1,2) {\footnotesize 1};
	\draw [fill] (1,2) circle [radius=0.1];
	\draw (0,0)--(1,1)--(1,2);
	\draw (1,1)--(2,0);
	\end{tikzpicture}};
\draw (1)--(6);
\draw (1)--(8);
\draw (1)--(10);
\draw (2)--(6);
\draw (2)--(7);
\draw (3)--(8);
\draw (3)--(9);
\draw (4)--(10);
\draw (4)--(11);
\draw (5)--(7);
\draw (5)--(9);
\draw (5)--(11);
\draw (6)--(12);
\draw (6)--(13);
\draw (7)--(12);
\draw (7)--(13);
\draw (8)--(12);
\draw (8)--(14);
\draw (9)--(12);
\draw (9)--(14);
\draw (10)--(13);
\draw (10)--(14);
\draw (11)--(13);
\draw (11)--(14);
\draw (12)--(15);
\draw (13)--(15);
\draw (14)--(15);
\end{tikzpicture}
\end{center}
\caption{The Tuffley poset $S([3])$.}
\end{figure}
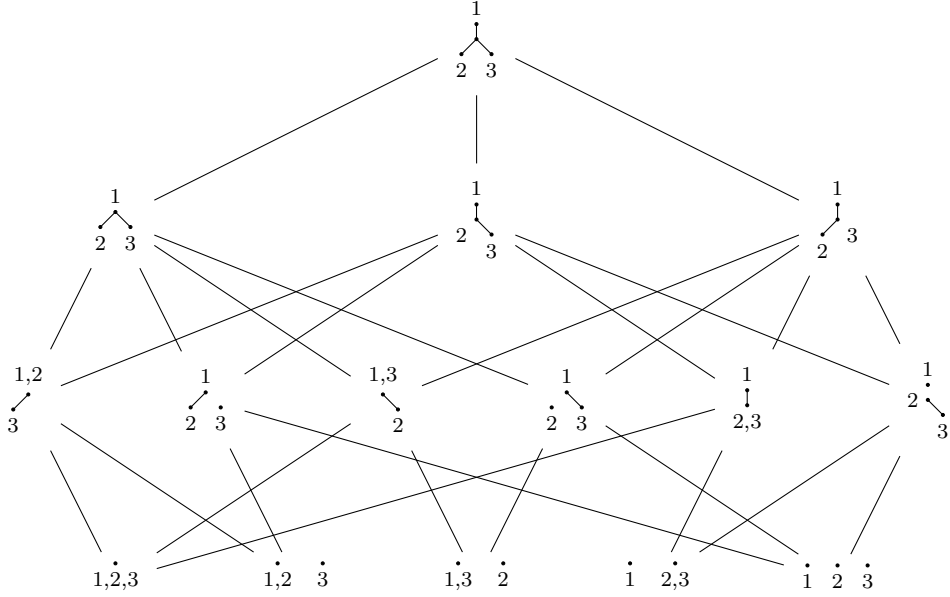


\subsection{Shellability for Regular CW Complexes}


We now turn to an overview of shellability for regular CW complexes. Following the convention of Bj{\"o}rner in \cite{bjorner}, we will call a finite, regular CW complex where each facet has dimension $d$ a \textit{$d$-CW-complex}. If $\mathcal{K}$ is a $d$-CW-complex and $f$ is a $d$-dimensional cell in $\mathcal{K}$, we will denote the $(d-1)$-CW-complex consisting of all proper faces of $f$ by $\partial f$. 

Gill, Linusson, Moulton and Steel proved in \cite{gill} that the CW decomposition given in \cite{oranges} is a regular CW complex.
The class of regular CW complexes behave like simplicial complexes in many ways; for example, regular CW complexes have face posets that determine their topology. In particular, given the face poset $P(\mathcal{K})$ of a regular CW complex $\mathcal{K}$, the order complex $\Delta(P(\mathcal{K})- \{\hat{0}\})$ is known as the \textit{barycentric subdivision of $\mathcal{K}$}; its geometric realization $||\Delta(P(\mathcal{K}) - \{\hat{0}\})||$  is homeomorphic to $\mathcal{K}$. Because of this, one is sometimes able to determine that a CW complex is shellable by examining the face poset of the complex. Thus the notion of shellability was extended from simplicial complexes to regular CW complexes by Bj{\"o}rner in  \cite{bjorner}. We state this generalization below. 

\begin{defn} A \textit{shelling} of the $d$-CW-complex $\mathcal{K}$ is an ordering $f_1, f_2, \ldots f_n$ of the $d$-dimensional cells of $\mathcal{K}$ where either $d=0$ or $d > 0$ and the following hold: 
\begin{enumerate}
\item $\partial f_i$ has a shelling
\item for $j=2, 3, \ldots n$, $\partial f_j \cap \left( \cup_{1 \leq i \leq j-1} \partial f_i \right)$ is a $(d-1)$-CW-complex 
\item for $j=2, 3, \ldots n$, $\partial f_j$ has a shelling in which the $(d-1)$-cells of $\partial f_j \cap \left(\cup_{1 \leq i \leq j-1} \partial f_i \right)$ come first. 
\end{enumerate}
If a shelling of $\mathcal{K}$ exists, $\mathcal{K}$ is said to be \textit{shellable}. 
\label{dcwshell}
\end{defn}

%

One technical but very convenient way of determining that a regular CW complex is shellable is by finding a recursive coatom ordering for the face poset of $\mathcal{K}$ augmented by a unique maximum, $\hat{1}$. This method was introduced by Bj{\"o}rner and Wachs in \cite{bw} as a recursive formulation of CL-shellability, a version of lexicographic shellability. 

\begin{defn} \label{rco} A bounded, graded poset $P$ admits a \textit{recursive coatom ordering} if there is an ordering $C_1, C_2, \ldots C_t$ of the coatoms of $P$ satisfying: 
\begin{enumerate}[label=(\roman*)]
\item For all $j=1, \ldots, t$, $[\hat{0}, C_j]$ admits a recursive coatom ordering with the property that for $j \neq 1$, the coatoms that come first in the ordering are those covered by some $C_k$ for $k < j$. 
\item For all $i < j$ and $y < C_i, C_j$, there exists a $k < j$ and an element $z$ such that $z \lessdot a_k$, $z \lessdot a_j$ and $y <  z$. 
\end{enumerate}

\end{defn}

The implication that the edge-product space is not shellable arises from the following result from \cite{bjorner}. Though Bj{\"o}rner originally presented this proposition in terms of CL-shellability, for our purposes, it will be convenient to present it in terms of recursive coatom orderings instead.

\begin{prop} [Proposition~4.2 \cite{bjorner}] A $d$-CW complex $\mathcal{K}$ is shellable if and only if its augmented face poset $\hat{P}(\mathcal{K})=P(\mathcal{K}) \cup \{\hat{1}\}$ admits a recursive coatom ordering. 
\label{shelliff}
\end{prop}


Note that shellability of a CW complex is not a topological property but is instead a property of the cell decomposition. Thus, though it was established in \cite{gill} that $\mathcal{E}(X) \cong ||\Delta(S(X))||$, it is possible that $||\Delta(S(X))||$ is shellable while $\mathcal{E}(X)$ is not. As mentioned earlier, examples of regular CW complexes that are not shellable but that have shellable barycentric subdivision are provided in \cite{vincewachs}, \cite{walker}, and \cite{vince}.

\subsection{NNI-tree Space}
We end this preliminary section with an overview of a graph known as NNI-tree space. This graph will be used to prove both that the edge-product space is gallery-connected (Proposition \ref{galconnectedthm}) and that the edge-product space is not shellable (Theorem \ref{s5norco}). 

NNI-tree space arises from a long-studied question in phylogenetics concerning optimal methods for comparing two different $X$-trees. Most methods define an operation on a tree and then define a metric between two trees (sometimes called a tree rearrangement metric) to be the minimum number of applications of this operation required to transform one tree into the other. For a general overview of the most common tree rearrangement metrics, see \cite{steelphylo}. The results in this paper focus specifically on a metric first introduced by Robinson in \cite{robinson71}, involving an operation on trivalent trees called crossovers.  Later, Waterman and Smith \cite{ws78} and others explored this operation which they describe as ``nearest neighbor interchange." We will also use this term to describe this operation: 

\begin{defn} Let $\alpha=(u, v)$ be any internal edge in $C_i \in \mathcal{C}_X$. The edge $\alpha$ divides subtrees $A$ and $B$ (each adjacent to vertex $u$) from subtrees $C$ and $D$ (each adjacent to vertex $v$). Swapping subtrees $B$ and $C$ yields a new tree $C_j \in \mathcal{C}_X$. Swapping subtrees $B$ and $D$ yields a new tree $C_k \in \mathcal{C}_X$. The operation of obtaining either $C_j$ or $C_k$ from $C_i$ is called a \textit{nearest neighbor interchange} (or \textit{NNI}) over edge $\alpha$. 
\label{crossover}
\end{defn} 

\begin{figure} 
\begin{center}
\begin{tikzpicture}[scale=0.8]
\node [below] at (0,0) {$A$};
\draw [fill] (0,0) circle [radius=0.05];
\node [below] at (2,0) {$B$};
\draw [fill] (2,0) circle [radius=0.05];
\node [right] at (1,1) {$u$};
\draw [fill] (1,1) circle [radius=0.05];
\node [right] at (2,2) {$v$};
\draw [fill] (2,2) circle [radius=0.05];
\node [right] at (3,3) {$D$};
\draw [fill] (3,3) circle [radius=0.05];
\draw [fill] (4,0) circle [radius=0.05];
\node [below] at (4,0) {$C$};
\draw (0,0)--(1,1)--(2,2) node [midway, above left ] {$\alpha$}--(3,3);
\draw (1,1)--(2,0);
\draw (2,2)--(4,0);
\end{tikzpicture}
\hskip 0.4cm
\begin{tikzpicture}[scale=0.8]
\node [below] at (0,0) {$A$};
\draw [fill] (0,0) circle [radius=0.05];
\node [below] at (2,0) {$C$};
\draw [fill] (2,0) circle [radius=0.05];
\node [right] at (1,1) {$u$};
\draw [fill] (1,1) circle [radius=0.05];
\node [right] at (2,2) {$v$};
\draw [fill] (2,2) circle [radius=0.05];
\node [right] at (3,3) {$D$};
\draw [fill] (3,3) circle [radius=0.05];
\draw [fill] (4,0) circle [radius=0.05];
\node [below] at (4,0) {$B$};
\draw (0,0)--(1,1) --(2,2) node [midway, above left ] {$\alpha$}--(3,3);
\draw (1,1)--(2,0);
\draw (2,2)--(4,0);
\end{tikzpicture}
\hskip 0.4cm
\begin{tikzpicture}[scale=0.8]
\node [below] at (0,0) {$A$};
\draw [fill] (0,0) circle [radius=0.05];
\node [below] at (2,0) {$D$};
\draw [fill] (2,0) circle [radius=0.05];
\node [right] at (1,1) {$u$};
\draw [fill] (1,1) circle [radius=0.05];
\node [right] at (2,2) {$v$};
\draw [fill] (2,2) circle [radius=0.05];
\node [right] at (3,3) {$B$};
\draw [fill] (3,3) circle [radius=0.05];
\draw [fill] (4,0) circle [radius=0.05];
\node [below] at (4,0) {$C$};
\draw (0,0)--(1,1) --(2,2) node [midway, above left ] {$\alpha$}--(3,3);
\draw (1,1)--(2,0);
\draw (2,2)--(4,0);
\end{tikzpicture}
\end{center}
\caption{Three maximal elements of the Tuffley poset, $S(X)$, all related by a NNI over edge $\alpha$.}
\label{crossoverpic} 
\end{figure}
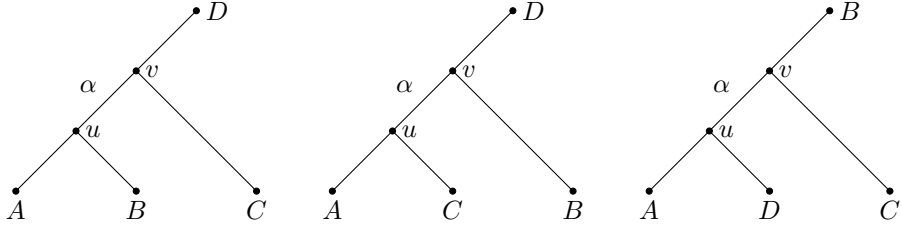

See Figure \ref{crossoverpic} for an example of three elements of $\mathcal{C}_X$ that are all related by a NNI over edge $\alpha$. 

\begin{defn} The \textit{NNI-tree space} is a graph with vertex set consisting of elements of $\mathcal{C}_X$ and an edge between $C_i$ and $C_j$ if and only if $C_i$ and $C_j$ are related by a NNI over some edge $\alpha$. When $|X|=n$, we denote this graph by $G_{NNI}(n)$. 
\end{defn} 


The upcoming proposition and corollary are important to proving results in Section \ref{noshell}, but we find it useful to introduce them here because they highlight fundamental properties of $S(X)$ and because they follow from Theorem \ref{tuffleyprops}; as such, we state them without proof. 
%

\begin{prop} Suppose $w \in S(X)$ is such that $w \lessdot C_i$ for some $C_i \in \mathcal{C}_X$. 
\begin{enumerate}[label=(\roman*)]
\item The element $w$ is covered by three elements if and only if $w$ is obtained by contracting an internal edge $\alpha$ in $C_i$. These three elements are related by an NNI over edge $\alpha$.
\item The element $w$ is covered by only $C_i$ if and only if $w$ is obtained by contracting a leaf edge in $C_i$. 
\end{enumerate}
\label{3coatomscoverw}
\end{prop}


\begin{cor} Suppose $C_i, C_j \in \mathcal{C}_X$.  Both $C_i$ and $C_j$ cover a common element $w$ in $S(X)$ if and only if $C_j$ can be obtained from $C_i$ by an NNI over some internal edge $\alpha$. 
\label{commonelement}
\end{cor} 
 
Since by Proposition \ref{3coatomscoverw}, $C_i$ and $C_j$ are related by a NNI over some edge $\alpha$ if and only if $C_i$ and $C_j$ cover a common element in $\hat{S}(X)$, NNI-tree space can be viewed as the facet-ridge graph of $\mathcal{E}(X)$. Recall that the facet-ridge graph of a $d$-CW-complex $\mathcal{K}$ has the $d$-dimensional facets of $\mathcal{K}$ as vertices and an edge between facets if they contain a common $(d-1)$-face of $\mathcal{K}$. Because NNI-tree space is the facet-ridge graph of $\mathcal{E}(X)$, we can use $G_{NNI}(n)$ to prove that $\mathcal{E}(X)$ is gallery connected. Shellable CW complexes are gallery-connected, though the converse is not necessarily true as we shall see in the case of $\mathcal{E}(X)$.
 
\begin{defn} A pure CW complex $\mathcal{K}$ is \textit{gallery-connected} if for any two $d$-dimensional facets $C_j$ and $C_k$, there exists a sequence of $d$-dimensional facets $C_j=C_0, C_1, C_2, \ldots, C_{r-1}, C_r=C_k$ where for $1 \leq i \leq r$, $C_{i-1}$ and $C_{i}$ contain a common $(d-1)$-dimensional face of $\mathcal{K}$. 
\label{galconnected}
\end{defn}

We use NNI-tree space to prove the following:

\begin{prop} 
\label{galconnectedthm}
The edge-product space $\mathcal{E}(X)$ is gallery-connected. 
\end{prop} 

\begin{proof}
Let  $|X|=n$ and $C_j$, $C_k$ be facets of $\mathcal{E}(X)$ so each corresponds to both a maximal element of $S(X)$ and a vertex in NNI-tree space. It is well-known that NNI-tree space is a connected graph; as a result, any path between the vertex corresponding to $C_j$ and the vertex corresponding to $C_k$ gives a sequence of vertices; each vertex in such a sequence corresponds to an $n$-dimensional face (facet) of $\mathcal{E}(X)$. If $v$ and $v'$ are adjacent vertices in $G_{NNI}(n)$, they correspond to maximal elements in $S(X)$ that cover a common element $u$ in $S(X)$ by definition of $G_{NNI}(n)$. The element $u$ corresponds to an $(n-1)$-dimensional face of $\mathcal{E}(X)$.  Thus this sequence of vertices gives the desired sequence of facets.
\end{proof}

\section{The poset $\hat{S}(X)$ does not admit a recursive coatom ordering}
\label{noshell}
This section serves to show that for $|X| \geq 5$, the edge-product space $\mathcal{E}(X)$ is not shellable. In particular, we will use Proposition \ref{shelliff} and show that $\hat{S}(X)=S(X) \cup \{\hat{0}, \hat{1}\}$ does not admit a recursive coatom ordering because for any total order on the coatoms of $\hat{S}(X)$, (ii) of Definition \ref{rco} fails. Note that recursive coatom orderings exist for $|X| \leq 4$.

\begin{figure}
\begin{center}
\begin{tikzpicture}[scale=0.75]
\node [] at (5, 15){\begin{tikzpicture}[scale=0.25] 
\node at (10,0) [below] {\tiny 1};
\draw [fill] (10,0)  circle [radius=0.05]; 
\node at (12,0) [below] {\tiny n-1};
\draw [fill] (12,0)  circle [radius=0.05]; 
\node at (14,0) [below] {\tiny n-2};
\draw [fill] (14,0)  circle [radius=0.05]; 
\node at (16,0) [below] {\tiny n-3};
\draw [fill] (16,0)  circle [radius=0.05]; 
\node at (11,1) [below] {};
\draw [fill] (11,1)  circle [radius=0.05]; 
\node at (15,6) [left] {\tiny n};
\draw [fill] (15,6)  circle [radius=0.05];
\node at (13,3) [below] {};
\draw [fill] (13,3)  circle [radius=0.05]; 
\node at (20,0) [below] {\tiny 2};
\draw [fill] (20,0)  circle [radius=0.05]; 
\draw [fill] (12,2)  circle [radius=0.05]; 
\draw [fill] (15,5)  circle [radius=0.05]; 
\node at (13,-3) {}; 
\draw [dashed] (13,3)--(15,5);
\draw (10,0)--(11,1);
\draw (12,0)--(11,1);
\draw (11,1)--(12,2);
\draw (14,0)--(12,2);
\draw (12,2)--(13,3);
\draw (16,0)--(13,3);
\draw (15,5)--(15,6);
\draw (20,0)--(15,5);
\end{tikzpicture}};
\node at (0, 10){\begin{tikzpicture}[scale=0.25] 
\node at (10,0) [below] {\tiny 1};
\draw [fill] (10,0)  circle [radius=0.05]; 
\node at (12,0) [below] {\tiny n-2};
\draw [fill] (12,0)  circle [radius=0.05]; 
\node at (14,0) [below] {\tiny n-1};
\draw [fill] (14,0)  circle [radius=0.05]; 
\node at (16,0) [below] {\tiny n-3};
\draw [fill] (16,0)  circle [radius=0.05]; 
\node at (11,1) [below] {};
\draw [fill] (11,1)  circle [radius=0.05]; 
\node at (15,6) [left] {\tiny n};
\draw [fill] (15,6)  circle [radius=0.05];
\node at (13,3) [below] {};
\draw [fill] (13,3)  circle [radius=0.05]; 
\node at (20,0) [below] {\tiny 2};
\draw [fill] (20,0)  circle [radius=0.05]; 
\draw [fill] (12,2)  circle [radius=0.05]; 
\draw [fill] (15,5)  circle [radius=0.05]; 
\node at (13.5,-3) {}; 
\draw [dashed] (13,3)--(15,5);
\draw (10,0)--(11,1);
\draw (12,0)--(11,1);
\draw (11,1)--(12,2);
\draw (14,0)--(12,2);
\draw (12,2)--(13,3);
\draw (16,0)--(13,3);
\draw (15,5)--(15,6);
\draw (20,0)--(15,5);
\end{tikzpicture}};
\node at (0, 5){\begin{tikzpicture}[scale=0.25] 
\node at (10,0) [below] {\tiny 1};
\draw [fill] (10,0)  circle [radius=0.05]; 
\node at (12,0) [below] {\tiny n-2};
\draw [fill] (12,0)  circle [radius=0.05]; 
\node at (14,0) [below] {\tiny n-3};
\draw [fill] (14,0)  circle [radius=0.05]; 
\node at (16,0) [below] {\tiny n-1};
\draw [fill] (16,0)  circle [radius=0.05]; 
\node at (11,1) [below] {};
\draw [fill] (11,1)  circle [radius=0.05]; 
\node at (15,6) [left] {\tiny n};
\draw [fill] (15,6)  circle [radius=0.05];
\node at (13,3) [below] {};
\draw [fill] (13,3)  circle [radius=0.05]; 
\node at (20,0) [below] {\tiny 2};
\draw [fill] (20,0)  circle [radius=0.05]; 
\draw [fill] (12,2)  circle [radius=0.05]; 
\draw [fill] (15,5)  circle [radius=0.05]; 
\node at (13.5,-3) {}; 
\draw [dashed] (13,3)--(15,5);
\draw (10,0)--(11,1);
\draw (12,0)--(11,1);
\draw (11,1)--(12,2);
\draw (14,0)--(12,2);
\draw (12,2)--(13,3);
\draw (16,0)--(13,3);
\draw (15,5)--(15,6);
\draw (20,0)--(15,5);
\end{tikzpicture}};
\node at (5, 0){\begin{tikzpicture}[scale=0.25] 
\node at (10,0) [below] {\tiny 1};
\draw [fill] (10,0)  circle [radius=0.05]; 
\node at (12,0) [below] {\tiny n-3};
\draw [fill] (12,0)  circle [radius=0.05]; 
\node at (14,0) [below] {\tiny n-2};
\draw [fill] (14,0)  circle [radius=0.05]; 
\node at (16,0) [below] {\tiny n-1};
\draw [fill] (16,0)  circle [radius=0.05]; 
\node at (11,1) [below] {};
\draw [fill] (11,1)  circle [radius=0.05]; 
\node at (15,6) [left] {\tiny n};
\draw [fill] (15,6)  circle [radius=0.05];
\node at (13,3) [below] {};
\draw [fill] (13,3)  circle [radius=0.05]; 
\node at (20,0) [below] {\tiny 2};
\draw [fill] (20,0)  circle [radius=0.05]; 
\draw [fill] (12,2)  circle [radius=0.05]; 
\draw [fill] (15,5)  circle [radius=0.05]; 
\node at (13.5,-3) {}; 
\draw [dashed] (13,3)--(15,5);
\draw (10,0)--(11,1);
\draw (12,0)--(11,1);
\draw (11,1)--(12,2);
\draw (14,0)--(12,2);
\draw (12,2)--(13,3);
\draw (16,0)--(13,3);
\draw (15,5)--(15,6);
\draw (20,0)--(15,5);
\end{tikzpicture}};
\node at (10, 5){\begin{tikzpicture}[scale=0.25] 
\node at (10,0) [below] {\tiny 1};
\draw [fill] (10,0)  circle [radius=0.05]; 
\node at (12,0) [below] {\tiny n-3};
\draw [fill] (12,0)  circle [radius=0.05]; 
\node at (14,0) [below] {\tiny n-1};
\draw [fill] (14,0)  circle [radius=0.05]; 
\node at (16,0) [below] {\tiny n-2};
\draw [fill] (16,0)  circle [radius=0.05]; 
\node at (11,1) [below] {};
\draw [fill] (11,1)  circle [radius=0.05]; 
\node at (15,6) [left] {\tiny n};
\draw [fill] (15,6)  circle [radius=0.05];
\node at (13,3) [below] {};
\draw [fill] (13,3)  circle [radius=0.05]; 
\node at (20,0) [below] {\tiny 2};
\draw [fill] (20,0)  circle [radius=0.05]; 
\draw [fill] (12,2)  circle [radius=0.05]; 
\draw [fill] (15,5)  circle [radius=0.05]; 
\node at (13.5,-3) {}; 
\draw [dashed] (13,3)--(15,5);
\draw (10,0)--(11,1);
\draw (12,0)--(11,1);
\draw (11,1)--(12,2);
\draw (14,0)--(12,2);
\draw (12,2)--(13,3);
\draw (16,0)--(13,3);
\draw (15,5)--(15,6);
\draw (20,0)--(15,5);
\end{tikzpicture}};
\node at (10, 10){\begin{tikzpicture}[scale=0.25] 
\node at (10,0) [below] {\tiny 1};
\draw [fill] (10,0)  circle [radius=0.05]; 
\node at (12,0) [below] {\tiny n-1};
\draw [fill] (12,0)  circle [radius=0.05]; 
\node at (14,0) [below] {\tiny n-3};
\draw [fill] (14,0)  circle [radius=0.05]; 
\node at (16,0) [below] {\tiny n-2};
\draw [fill] (16,0)  circle [radius=0.05]; 
\node at (11,1) [below] {};
\draw [fill] (11,1)  circle [radius=0.05]; 
\node at (15,6) [left] {\tiny n};
\draw [fill] (15,6)  circle [radius=0.05];
\node at (13,3) [below] {};
\draw [fill] (13,3)  circle [radius=0.05]; 
\node at (20,0) [below] {\tiny 2};
\draw [fill] (20,0)  circle [radius=0.05]; 
\draw [fill] (12,2)  circle [radius=0.05]; 
\draw [fill] (15,5)  circle [radius=0.05]; 
\node at (13.5,-3) {}; 
\draw [dashed] (13,3)--(15,5);
\draw (10,0)--(11,1);
\draw (12,0)--(11,1);
\draw (11,1)--(12,2);
\draw (14,0)--(12,2);
\draw (12,2)--(13,3);
\draw (16,0)--(13,3);
\draw (15,5)--(15,6);
\draw (20,0)--(15,5);
\end{tikzpicture}};
\draw (5, 13.9)--(0,11.6);
\draw (0, 8.9)--(0, 6.6);
\draw (0, 3.9)--(5, 1.6);
\draw (5, 1.6)--(10, 3.9);
\draw (10, 6.6)--(10, 8.9);
\draw (10, 11.6)--(5, 13.9);
\end{tikzpicture}
\end{center}
\caption{The elements in the set $K$ form a nontrivial cycle in $G_{NNI}(n)$ for $n \geq 5$.}
\label{3cycleex}
\end{figure}
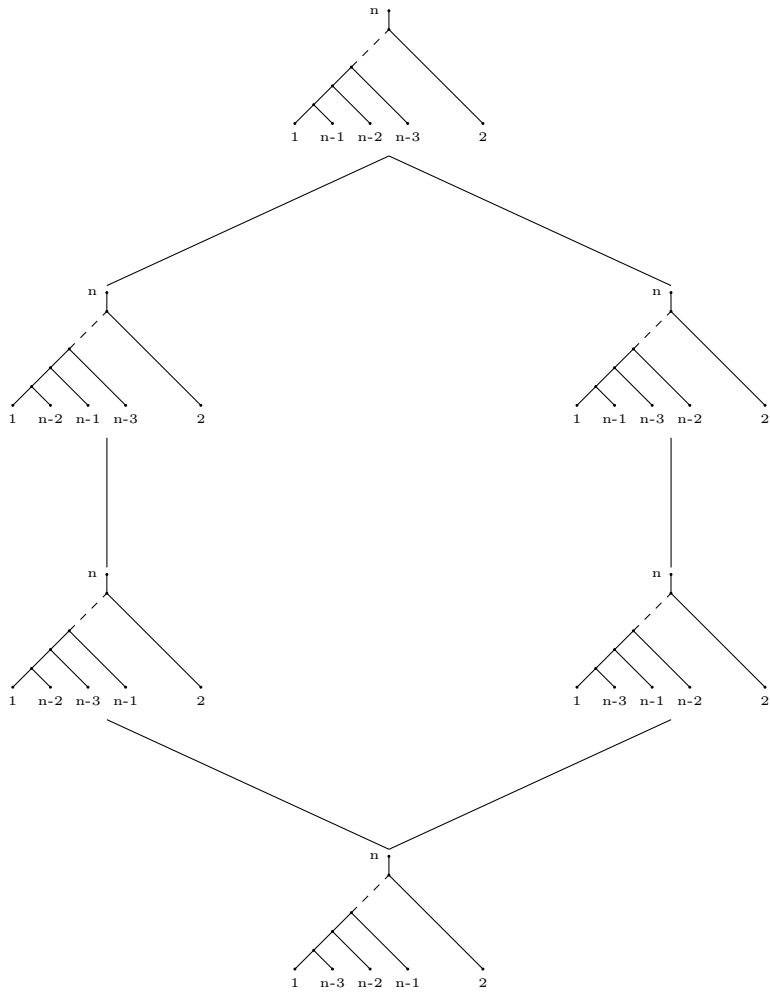

The general idea for this proof is as follows. First, we will simplify matters by taking $X=[n]$. We show that for $n \geq 5$ and for any ordering $\Omega$ on the coatoms of $\hat{S}([n])$, one of the elements in Figure \ref{3cycleex} comes last, say $C_i$. The two elements adjacent to $C_i$ in Figure \ref{3cycleex}, which we will call $C_j$ and $C_k$, are both above an element $F_{ijk}$ (defined in Definition \ref{problemdef}). Supposing without loss of generality that $C_j$ comes after $C_k$ in $\Omega$, we prove that there is no coatom that comes before $C_j$ in $\Omega$ that covers a common element $z$ with $C_j$, where $z$ is above $F_{ijk}$ in $\hat{S}([n])$. 

Lemma \ref{tfstructure} and Corollary \ref{tfstructurecor} are instrumental in proving this fact; both results follow from the idea that we can associate any element that covers $F_{ijk}$ in $\hat{S}([n])$ with a vertex of a tree $T(F_{ijk})$ and we can associate any coatom above $F_{ijk}$ in $\hat{S}([n])$ with an edge of the tree $T(F_{ijk})$. The tree $T(F_{ijk})$ has the property that two coatoms above $F_{ijk}$ are associated with adjacent edges in $T(F_{ijk})$ if and only if they cover a common element and this common element covers $F_{ijk}$ in $\hat{S}([n])$.  

%
%

The remainder of this section introduces concepts and notation that are used to prove Theorem \ref{s5norco}. We provide a complete example of the argument against the existence of a recursive coatom ordering for $\hat{S}([5])$ in Example \ref{s5ex}.

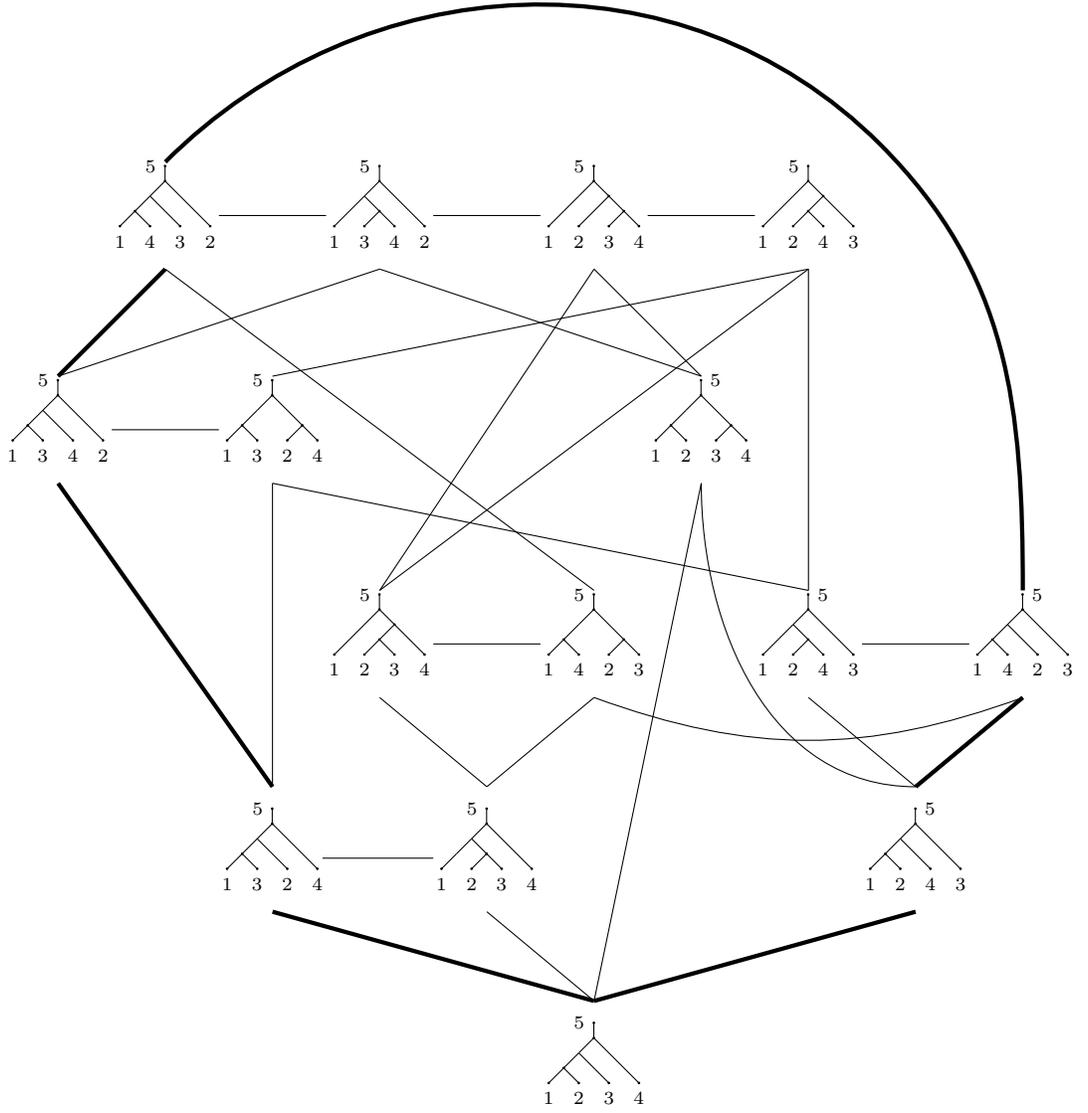
\begin{figure}
\begin{center}
\begin{tikzpicture}[scale=0.95]
\node at (4.5,3) {\begin{tikzpicture}[scale=0.2] \node at (10,0) [below] {\scriptsize 1};
\draw [fill] (10,0)  circle [radius=0.05]; 
\node at (12,0) [below] {\scriptsize 2};
\draw [fill] (12,0)  circle [radius=0.05]; 
\node at (14,0) [below] {\scriptsize 3};
\draw [fill] (14,0)  circle [radius=0.05]; 
\node at (16,0) [below] {\scriptsize 4};
\draw [fill] (16,0)  circle [radius=0.05]; 
\node at (11,1) [below] {};
\draw [fill] (11,1)  circle [radius=0.05]; 
\node at (13,4) [left] {\scriptsize 5};
\draw [fill] (13,4)  circle [radius=0.05];
\node at (13,3) [below] {};
\draw [fill] (13,3)  circle [radius=0.05]; 
\node at (13.5,-3) {}; 
\draw (10,0)--(11,1);
\draw (12,0)--(11,1);
\draw (11,1)--(12,2);
\draw (14,0)--(12,2);
\draw (12,2)--(13,3);
\draw (16,0)--(13,3);
\draw (13,3)--(13,4);
\end{tikzpicture}};
\node at (0,6) {\begin{tikzpicture}[scale=0.2] \node at (10,0) [below] {\scriptsize 1};
\draw [fill] (10,0)  circle [radius=0.05]; 
\node at (12,0) [below] {\scriptsize 3};
\draw [fill] (12,0)  circle [radius=0.05]; 
\node at (14,0) [below] {\scriptsize 2};
\draw [fill] (14,0)  circle [radius=0.05]; 
\node at (16,0) [below] {\scriptsize 4};
\draw [fill] (16,0)  circle [radius=0.05]; 
\node at (11,1) [below] {};
\draw [fill] (11,1)  circle [radius=0.05]; 
\node at (13,4) [left] {\scriptsize 5};
\draw [fill] (13,4)  circle [radius=0.05];
\node at (13,3) [below] {};
\draw [fill] (13,3)  circle [radius=0.05]; 
\node at (13.5,-3) {}; 
\draw (10,0)--(11,1);
\draw (12,0)--(11,1);
\draw (11,1)--(12,2);
\draw (14,0)--(12,2);
\draw (12,2)--(13,3);
\draw (16,0)--(13,3);
\draw (13,3)--(13,4);
\end{tikzpicture}};;
\node at (3, 6){\begin{tikzpicture}[scale=0.2] \node at (10,0) [below] {\scriptsize 1};
\draw [fill] (10,0)  circle [radius=0.05]; 
\node at (12,0) [below] {\scriptsize 2};
\draw [fill] (12,0)  circle [radius=0.05]; 
\node at (14,0) [below] {\scriptsize 3};
\draw [fill] (14,0)  circle [radius=0.05]; 
\node at (16,0) [below] {\scriptsize 4};
\draw [fill] (16,0)  circle [radius=0.05]; 
\node at (11,1) [below] {};
\node at (13,4) [left] {\scriptsize 5};
\draw [fill] (13,4)  circle [radius=0.05];
\node at (13,3) [below] {};
\draw [fill] (13,3)  circle [radius=0.05]; 
\draw [fill] (13,1)  circle [radius=0.05];
\node at (13.5,-3) {}; 
\draw (10,0)--(11,1);
\draw (12,0)--(13,1);
\draw (11,1)--(12,2);
\draw (14,0)--(12,2);
\draw (12,2)--(13,3);
\draw (16,0)--(13,3);
\draw (13,3)--(13,4);
\end{tikzpicture}};
\node at (6,12)  {\begin{tikzpicture}[scale=0.2] \node at (10,0) [below] {\scriptsize 1};
\draw [fill] (10,0)  circle [radius=0.05]; 
\node at (12,0) [below] {\scriptsize 2};
\draw [fill] (12,0)  circle [radius=0.05]; 
\node at (14,0) [below] {\scriptsize 3};
\draw [fill] (14,0)  circle [radius=0.05]; 
\node at (16,0) [below] {\scriptsize 4};
\draw [fill] (16,0)  circle [radius=0.05]; 
\node at (11,1) [below] {};
\draw [fill] (11,1)  circle [radius=0.05]; 
\node at (13,4) [right] {\scriptsize 5};
\draw [fill] (13,4)  circle [radius=0.05];
\node at (13,3) [below] {};
\draw [fill] (13,3)  circle [radius=0.05]; 
\draw [fill] (15,1) circle [radius=0.05];
\node at (13.5,-3) {}; 
\draw (10,0)--(11,1);
\draw (12,0)--(11,1);
\draw (11,1)--(12,2);
\draw (14,0)--(15,1);
\draw (12,2)--(13,3);
\draw (16,0)--(13,3);
\draw (13,3)--(13,4);
\end{tikzpicture}};;
\node at (9, 6) {\begin{tikzpicture}[scale=0.2] \node at (10,0) [below] {\scriptsize 1};
\draw [fill] (10,0)  circle [radius=0.05]; 
\node at (12,0) [below] {\scriptsize 2};
\draw [fill] (12,0)  circle [radius=0.05]; 
\node at (14,0) [below] {\scriptsize 4};
\draw [fill] (14,0)  circle [radius=0.05]; 
\node at (16,0) [below] {\scriptsize 3};
\draw [fill] (16,0)  circle [radius=0.05]; 
\node at (11,1) [below] {};
\draw [fill] (11,1)  circle [radius=0.05]; 
\node at (13,4) [right] {\scriptsize 5};
\draw [fill] (13,4)  circle [radius=0.05];
\node at (13,3) [below] {};
\draw [fill] (13,3)  circle [radius=0.05]; 
\node at (13.5,-3) {}; 
\draw (10,0)--(11,1);
\draw (12,0)--(11,1);
\draw (11,1)--(12,2);
\draw (14,0)--(12,2);
\draw (12,2)--(13,3);
\draw (16,0)--(13,3);
\draw (13,3)--(13,4);
\end{tikzpicture}};;
\node at (-3,12) {\begin{tikzpicture}[scale=0.2] \node at (10,0) [below] {\scriptsize 1};
\draw [fill] (10,0)  circle [radius=0.05]; 
\node at (12,0) [below] {\scriptsize 3};
\draw [fill] (12,0)  circle [radius=0.05]; 
\node at (14,0) [below] {\scriptsize 4};
\draw [fill] (14,0)  circle [radius=0.05]; 
\node at (16,0) [below] {\scriptsize 2};
\draw [fill] (16,0)  circle [radius=0.05]; 
\node at (11,1) [below] {};
\draw [fill] (11,1)  circle [radius=0.05]; 
\node at (13,4) [left] {\scriptsize 5};
\draw [fill] (13,4)  circle [radius=0.05];
\node at (13,3) [below] {};
\draw [fill] (13,3)  circle [radius=0.05]; 
\node at (13.5,-3) {}; 
\draw (10,0)--(11,1);
\draw (12,0)--(11,1);
\draw (11,1)--(12,2);
\draw (14,0)--(12,2);
\draw (12,2)--(13,3);
\draw (16,0)--(13,3);
\draw (13,3)--(13,4);
\end{tikzpicture}};;
\node at (0, 12) {\begin{tikzpicture}[scale=0.2] \node at (10,0) [below] {\scriptsize 1};
\draw [fill] (10,0)  circle [radius=0.05]; 
\node at (12,0) [below] {\scriptsize 3};
\draw [fill] (12,0)  circle [radius=0.05]; 
\node at (14,0) [below] {\scriptsize 2};
\draw [fill] (14,0)  circle [radius=0.05]; 
\node at (16,0) [below] {\scriptsize 4};
\draw [fill] (16,0)  circle [radius=0.05]; 
\node at (11,1) [below] {};
\draw [fill] (11,1)  circle [radius=0.05]; 
\node at (13,4) [left] {\scriptsize 5};
\draw [fill] (13,4)  circle [radius=0.05];
\node at (13,3) [below] {};
\draw [fill] (13,3)  circle [radius=0.05]; 
\draw [fill] (15,1) circle [radius=0.05];
\node at (13.5,-3) {}; 
\draw (10,0)--(11,1);
\draw (12,0)--(11,1);
\draw (11,1)--(12,2);
\draw (14,0)--(15,1);
\draw (12,2)--(13,3);
\draw (16,0)--(13,3);
\draw (13,3)--(13,4);
\end{tikzpicture}};;
\node at (4.5, 9){\begin{tikzpicture}[scale=0.2] \node at (10,0) [below] {\scriptsize 1};
\draw [fill] (10,0)  circle [radius=0.05]; 
\node at (12,0) [below] {\scriptsize 4};
\draw [fill] (12,0)  circle [radius=0.05]; 
\node at (14,0) [below] {\scriptsize 2};
\draw [fill] (14,0)  circle [radius=0.05]; 
\node at (16,0) [below] {\scriptsize 3};
\draw [fill] (16,0)  circle [radius=0.05]; 
\node at (11,1) [below] {};
\draw [fill] (11,1)  circle [radius=0.05]; 
\node at (13,4) [left] {\scriptsize 5};
\draw [fill] (13,4)  circle [radius=0.05];
\node at (13,3) [below] {};
\draw [fill] (13,3)  circle [radius=0.05]; 
\draw [fill] (15,1) circle [radius=0.05];
\node at (13.5,-3) {}; 
\draw (10,0)--(11,1);
\draw (12,0)--(11,1);
\draw (11,1)--(12,2);
\draw (14,0)--(15,1);
\draw (12,2)--(13,3);
\draw (16,0)--(13,3);
\draw (13,3)--(13,4);
\end{tikzpicture}}; 
\node at (1.5, 9){\begin{tikzpicture}[scale=0.2] \node at (10,0) [below] {\scriptsize 1};
\draw [fill] (10,0)  circle [radius=0.05]; 
\node at (12,0) [below] {\scriptsize 2};
\draw [fill] (12,0)  circle [radius=0.05]; 
\node at (14,0) [below] {\scriptsize 3};
\draw [fill] (14,0)  circle [radius=0.05]; 
\node at (16,0) [below] {\scriptsize 4};
\draw [fill] (16,0)  circle [radius=0.05]; 
\node at (11,1) [below] {};
\draw [fill] (14,2)  circle [radius=0.05]; 
\node at (13,4) [left] {\scriptsize 5};
\draw [fill] (13,4)  circle [radius=0.05];
\node at (13,3) [below] {};
\draw [fill] (13,3)  circle [radius=0.05]; 
\draw [fill] (13,1) circle [radius=0.05];
\node at (13.5,-3) {}; 
\draw (10,0)--(11,1);
\draw (12,0)--(14,2);
\draw (11,1)--(12,2);
\draw (14,0)--(13,1);
\draw (12,2)--(13,3);
\draw (16,0)--(13,3);
\draw (13,3)--(13,4);
\end{tikzpicture}};
\node at (7.5, 9){\begin{tikzpicture}[scale=0.2] \node at (10,0) [below] {\scriptsize 1};
\draw [fill] (10,0)  circle [radius=0.05]; 
\node at (12,0) [below] {\scriptsize 2};
\draw [fill] (12,0)  circle [radius=0.05]; 
\node at (14,0) [below] {\scriptsize 4};
\draw [fill] (14,0)  circle [radius=0.05]; 
\node at (16,0) [below] {\scriptsize 3};
\draw [fill] (16,0)  circle [radius=0.05]; 
\node at (11,1) [below] {};
\draw [fill] (13,1)  circle [radius=0.05]; 
\node at (13,4) [right] {\scriptsize 5};
\draw [fill] (13,4)  circle [radius=0.05];
\node at (13,3) [below] {};
\draw [fill] (13,3)  circle [radius=0.05]; 
\node at (13.5,-3) {}; 
\draw (10,0)--(11,1);
\draw (12,0)--(13,1);
\draw (11,1)--(12,2);
\draw (14,0)--(12,2);
\draw (12,2)--(13,3);
\draw (16,0)--(13,3);
\draw (13,3)--(13,4);
\end{tikzpicture}};;
\node at (10.5, 9) {\begin{tikzpicture}[scale=0.2] \node at (10,0) [below] {\scriptsize 1};
\draw [fill] (10,0)  circle [radius=0.05]; 
\node at (12,0) [below] {\scriptsize 4};
\draw [fill] (12,0)  circle [radius=0.05]; 
\node at (14,0) [below] {\scriptsize 2};
\draw [fill] (14,0)  circle [radius=0.05]; 
\node at (16,0) [below] {\scriptsize 3};
\draw [fill] (16,0)  circle [radius=0.05]; 
\node at (11,1) [below] {};
\draw [fill] (11,1)  circle [radius=0.05]; 
\node at (13,4) [right] {\scriptsize 5};
\draw [fill] (13,4)  circle [radius=0.05];
\node at (13,3) [below] {};
\draw [fill] (13,3)  circle [radius=0.05]; 
\node at (13.5,-3) {}; 
\draw (10,0)--(11,1);
\draw (12,0)--(11,1);
\draw (11,1)--(12,2);
\draw (14,0)--(12,2);
\draw (12,2)--(13,3);
\draw (16,0)--(13,3);
\draw (13,3)--(13,4);
\end{tikzpicture}};;;
\node at (4.5, 15) {\begin{tikzpicture}[scale=0.2] \node at (10,0) [below] {\scriptsize 1};
\draw [fill] (10,0)  circle [radius=0.05]; 
\node at (12,0) [below] {\scriptsize 2};
\draw [fill] (12,0)  circle [radius=0.05]; 
\node at (14,0) [below] {\scriptsize 3};
\draw [fill] (14,0)  circle [radius=0.05]; 
\node at (16,0) [below] {\scriptsize 4};
\draw [fill] (16,0)  circle [radius=0.05]; 
\node at (11,1) [below] {};
\draw [fill] (14,2)  circle [radius=0.05]; 
\node at (13,4) [left] {\scriptsize 5};
\draw [fill] (13,4)  circle [radius=0.05];
\node at (13,3) [below] {};
\draw [fill] (13,3)  circle [radius=0.05]; 
\draw [fill] (15,1) circle [radius=0.05];
\node at (13.5,-3) {}; 
\draw (10,0)--(11,1);
\draw (12,0)--(14,2);
\draw (11,1)--(12,2);
\draw (14,0)--(15,1);
\draw (12,2)--(13,3);
\draw (16,0)--(13,3);
\draw (13,3)--(13,4);
\end{tikzpicture}};
\node at (1.5, 15) {\begin{tikzpicture}[scale=0.2] \node at (10,0) [below] {\scriptsize 1};
\draw [fill] (10,0)  circle [radius=0.05]; 
\node at (12,0) [below] {\scriptsize 3};
\draw [fill] (12,0)  circle [radius=0.05]; 
\node at (14,0) [below] {\scriptsize 4};
\draw [fill] (14,0)  circle [radius=0.05]; 
\node at (16,0) [below] {\scriptsize 2};
\draw [fill] (16,0)  circle [radius=0.05]; 
\node at (11,1) [below] {};
\draw [fill] (13,1)  circle [radius=0.05]; 
\node at (13,4) [left] {\scriptsize 5};
\draw [fill] (13,4)  circle [radius=0.05];
\node at (13,3) [below] {};
\draw [fill] (13,3)  circle [radius=0.05]; 
\node at (13.5,-3) {}; 
\draw (10,0)--(11,1);
\draw (12,0)--(13,1);
\draw (11,1)--(12,2);
\draw (14,0)--(12,2);
\draw (12,2)--(13,3);
\draw (16,0)--(13,3);
\draw (13,3)--(13,4);
\end{tikzpicture}};
\node at (-1.5, 15)  {\begin{tikzpicture}[scale=0.2] \node at (10,0) [below] {\scriptsize 1};
\draw [fill] (10,0)  circle [radius=0.05]; 
\node at (12,0) [below] {\scriptsize 4};
\draw [fill] (12,0)  circle [radius=0.05]; 
\node at (14,0) [below] {\scriptsize 3};
\draw [fill] (14,0)  circle [radius=0.05]; 
\node at (16,0) [below] {\scriptsize 2};
\draw [fill] (16,0)  circle [radius=0.05]; 
\node at (11,1) [below] {};
\draw [fill] (11,1)  circle [radius=0.05]; 
\node at (13,4) [left] {\scriptsize 5};
\draw [fill] (13,4)  circle [radius=0.05];
\node at (13,3) [below] {};
\draw [fill] (13,3)  circle [radius=0.05]; 
\node at (13.5,-3) {}; 
\draw (10,0)--(11,1);
\draw (12,0)--(11,1);
\draw (11,1)--(12,2);
\draw (14,0)--(12,2);
\draw (12,2)--(13,3);
\draw (16,0)--(13,3);
\draw (13,3)--(13,4);
\end{tikzpicture}};;
\node at (7.5, 15) {\begin{tikzpicture}[scale=0.2] \node at (10,0) [below] {\scriptsize 1};
\draw [fill] (10,0)  circle [radius=0.05]; 
\node at (12,0) [below] {\scriptsize 2};
\draw [fill] (12,0)  circle [radius=0.05]; 
\node at (14,0) [below] {\scriptsize 4};
\draw [fill] (14,0)  circle [radius=0.05]; 
\node at (16,0) [below] {\scriptsize 3};
\draw [fill] (16,0)  circle [radius=0.05]; 
\node at (11,1) [below] {};
\draw [fill] (14,2)  circle [radius=0.05]; 
\node at (13,4) [left] {\scriptsize 5};
\draw [fill] (13,4)  circle [radius=0.05];
\node at (13,3) [below] {};
\draw [fill] (13,3)  circle [radius=0.05]; 
\draw [fill] (13,1) circle [radius=0.05];
\node at (13.5,-3) {}; 
\draw (10,0)--(11,1);
\draw (12,0)--(14,2);
\draw (11,1)--(12,2);
\draw (14,0)--(13,1);
\draw (12,2)--(13,3);
\draw (16,0)--(13,3);
\draw (13,3)--(13,4);
\end{tikzpicture}};
\draw [ultra thick] (4.5, 4)--(0,5.25);
\draw (4.5, 4)--(3,5.25);
\draw (4.5, 4)--(6, 11.25);
\draw [ultra thick](4.5, 4)--(9,5.25);
\draw [ultra thick](0, 7)--(-3, 11.25); 
\draw (0,7)--(0, 11.25); 
\draw (3,7)--(4.5, 8.25);
\draw (3,7)--(1.5, 8.25);
\draw (9, 7)--(7.5, 8.25);
\draw [ultra thick](9,7)--(10.5, 8.25);
\draw (2.25,6)--(0.7,6); 
\draw (8.25, 9)--(9.75, 9); 
\draw (2.25, 9)--(3.75, 9);
\draw (-2.25, 12)--(-0.75, 12);
\draw (-0.75, 15)--(0.75, 15); 
\draw (2.25, 15)--(3.75, 15); 
\draw (5.25, 15)--(6.75, 15); 
\draw (7.5, 9.75)--(0, 11.25); 
\draw (7.5, 9.75)--(7.5, 14.25); 
\draw (0, 12.75)--(7.5, 14.25); 
\draw [ultra thick] (-3, 12.75)--(-1.5, 14.25);
\draw (-3, 12.75)--(1.5, 14.25); 
\draw (1.5, 9.75)--(4.5, 14.25);
\draw (1.5, 9.75)--(7.5, 14.25);
\draw (4.5, 8.25) to [out=340, in=200] (10.5, 8.25);
\draw (4.5, 9.75)--(-1.5,14.25);
\draw (6, 12.75)--(4.5, 14.25);
\draw (6, 12.75)--(1.5, 14.25); 
\draw (9, 7) to [out=180, in=270] (6, 11.25);
\draw [ultra thick] (10.5, 9.75) to [out=90, in=315] (8.5, 16) to [out=135, in=45] (-1.5,15.75);
\end{tikzpicture}
\end{center}
\caption{The graph $G_{NNI}(5)$. The bold edges mark a nontrivial cycle of $G_{NNI}(5)$.}
\label{H5}
\end{figure}

We will show that for $n \geq 5$, a recursive coatom ordering for $\hat{S}([n])$ cannot exist because of the existence of certain triples of coatoms of $\hat{S}([n])$, which we describe next.

%

\begin{defn}
A \textit{critical triplet} is a triple $\{C_i, C_j, C_k\}$, where $C_i$, $C_j$, and $C_k \in \mathcal{C}_{[n]}$, satisfying the following: 

\begin{enumerate}
\item $C_i$ contains a leaf labeled $x$ such that the edge adjacent to this leaf, $e_x$ is adjacent to two internal edges, say $e_{x_1}$ and $e_{x_2}$, 
\item $C_i$ and $C_j$ cover $F_j$, the element obtained from $C_i$ by contracting $e_{x_1}$, and
\item $C_i$ and $C_k$ cover $F_k$, the element obtained from $C_i$ by contracting $e_{x_2}$.
\end{enumerate}
\label{tripletcond}

\end{defn}


\begin{prop} If $\{C_i, C_j, C_k\}$ is a critical triplet, then $C_j$ and $C_k$ are not adjacent in $G_{NNI}(n)$.  
\label{nonadjacent}
\end{prop}
\begin{proof} Note that for a critical triplet $\{C_i, C_j, C_k\}$, Proposition \ref{3coatomscoverw} implies that $C_i$ and $C_j$ (respectively $C_k$) are related by an NNI over $e_{x_1}$ (respectively $e_{x_2}$). Since $e_{x_1} \neq e_{x_2}$, it is straightforward to verify that $C_k$ can be obtained from $C_j$ by no fewer than two NNIs. Thus $C_j$ and $C_k$ do not cover a common element by Corollary \ref{commonelement}. 
\end{proof}

Given any critical triplet $\{C_i, C_j, C_k\}$, the intersection of the principal order ideals $I(C_j)$ and $I(C_k)$ contains an element that is crucial to our proof of Theorem \ref{s5norco}. We formally define this element next.

\begin{defn} 
\label{problemdef} 
Let $\{C_i, C_j, C_k\}$ be a critical triplet. The element $F_{ijk}$ is the element of $\hat{S}([n])$ obtained from $C_j$ by contracting $e_{x_1}$ then deleting $e_x$. Equivalently, $F_{ijk}$ can be obtained from $C_k$ by contracting $e_{x_2}$ then deleting $e_x$. 
\end{defn}

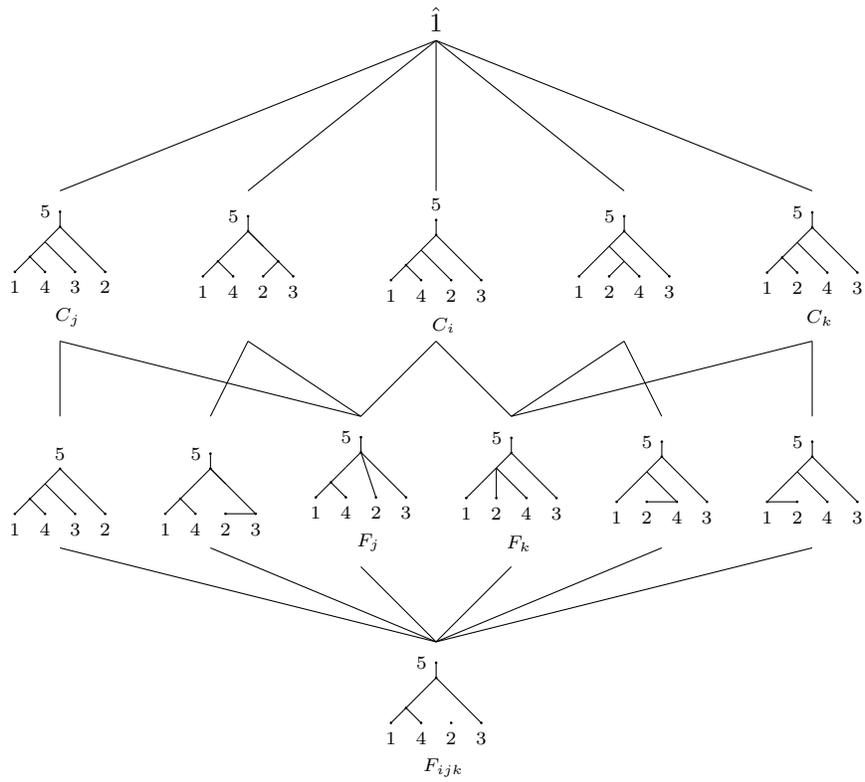
\begin{figure}
\begin{center}
\begin{tikzpicture}
\node at (0,6){
\begin{tikzpicture}[scale=0.2] \node at (10,0) [below] {\scriptsize 1};
\draw [fill] (10,0)  circle [radius=0.05]; 
\node at (12,0) [below] {\scriptsize 4};
\draw [fill] (12,0)  circle [radius=0.05]; 
\node at (14,0) [below] {\scriptsize 3};
\draw [fill] (14,0)  circle [radius=0.05]; 
\node at (16,0) [below] {\scriptsize 2};
\draw [fill] (16,0)  circle [radius=0.05]; 
\node at (11,1) [below] {};
\draw [fill] (11,1)  circle [radius=0.05]; 
\node at (13,4) [left] {\scriptsize 5};
\draw [fill] (13,4)  circle [radius=0.05];
\node at (13,3) [below] {};
\draw [fill] (13,3)  circle [radius=0.05]; 
\node at (13.5,-3) {\scriptsize $C_j$}; 
\draw (10,0)--(11,1);
\draw (12,0)--(11,1);
\draw (11,1)--(12,2);
\draw (14,0)--(12,2);
\draw (12,2)--(13,3);
\draw (16,0)--(13,3);
\draw (13,3)--(13,4);
\end{tikzpicture}};
\node at (10,6){
\begin{tikzpicture}[scale=0.2] \node at (10,0) [below] {\scriptsize 1};
\draw [fill] (10,0)  circle [radius=0.05]; 
\node at (12,0) [below] {\scriptsize 2};
\draw [fill] (12,0)  circle [radius=0.05]; 
\node at (14,0) [below] {\scriptsize 4};
\draw [fill] (14,0)  circle [radius=0.05]; 
\node at (16,0) [below] {\scriptsize 3};
\draw [fill] (16,0)  circle [radius=0.05]; 
\node at (11,1) [below] {};
\draw [fill] (11,1)  circle [radius=0.05]; 
\node at (13,4) [left] {\scriptsize 5};
\draw [fill] (13,4)  circle [radius=0.05];
\node at (13,3) [below] {};
\draw [fill] (13,3)  circle [radius=0.05]; 
\node at (13.5,-3) {\scriptsize $C_k$}; 
\draw (10,0)--(11,1);
\draw (12,0)--(11,1);
\draw (11,1)--(12,2);
\draw (14,0)--(12,2);
\draw (12,2)--(13,3);
\draw (16,0)--(13,3);
\draw (13,3)--(13,4);
\end{tikzpicture}};
\node at (2.5,6){\begin{tikzpicture}[scale=0.2] \node at (10,0) [below] {\scriptsize 1};
\draw [fill] (10,0)  circle [radius=0.05]; 
\node at (12,0) [below] {\scriptsize 4};
\draw [fill] (12,0)  circle [radius=0.05]; 
\node at (14,0) [below] {\scriptsize 2};
\draw [fill] (14,0)  circle [radius=0.05]; 
\node at (16,0) [below] {\scriptsize 3};
\draw [fill] (16,0)  circle [radius=0.05]; 
\node at (11,1) [below] {};
\draw [fill] (11,1)  circle [radius=0.05]; 
\node at (13,4) [left] {\scriptsize 5};
\draw [fill] (13,4)  circle [radius=0.05];
\node at (13,3) [below] {};
\draw [fill] (13,3)  circle [radius=0.05]; 
\node at (13.5,-3) {}; 
\draw [fill] (15,1) circle [radius=0.05];
\draw (10,0)--(11,1);
\draw (12,0)--(11,1);
\draw (11,1)--(12,2);
\draw (14,0)--(15,1);
\draw (12,2)--(13,3);
\draw (13,3)--(14,2);
\draw (16,0)--(13,3);
\draw (13,3)--(13,4);
\end{tikzpicture}};
\node at (5, 6){\begin{tikzpicture}[scale=0.2] \node at (10,0) [below] {\scriptsize 1};
\draw [fill] (10,0)  circle [radius=0.05]; 
\node at (12,0) [below] {\scriptsize 4};
\draw [fill] (12,0)  circle [radius=0.05]; 
\node at (14,0) [below] {\scriptsize 2};
\draw [fill] (14,0)  circle [radius=0.05]; 
\node at (16,0) [below] {\scriptsize 3};
\draw [fill] (16,0)  circle [radius=0.05]; 
\node at (11,1) [below] {};
\draw [fill] (11,1)  circle [radius=0.05]; 
\node at (13,4) [above] {\scriptsize 5};
\draw [fill] (13,4)  circle [radius=0.05];
\node at (13,3) [below] {};
\draw [fill] (13,3)  circle [radius=0.05]; 
\node at (13.5,-3) {\scriptsize $C_i$}; 
\draw (10,0)--(11,1);
\draw (12,0)--(11,1);
\draw (11,1)--(12,2);
\draw (14,0)--(12,2);
\draw (12,2)--(13,3);
\draw (16,0)--(13,3);
\draw (13,3)--(13,4);
\end{tikzpicture}};
\node at (7.5, 6) {\begin{tikzpicture}[scale=0.2] \node at (10,0) [below] {\scriptsize 1};
\draw [fill] (10,0)  circle [radius=0.05]; 
\node at (12,0) [below] {\scriptsize 2};
\draw [fill] (12,0)  circle [radius=0.05]; 
\node at (14,0) [below] {\scriptsize 4};
\draw [fill] (14,0)  circle [radius=0.05]; 
\node at (16,0) [below] {\scriptsize 3};
\draw [fill] (16,0)  circle [radius=0.05]; 
\node at (11,1) [below] {};
\node at (13,4) [left] {\scriptsize 5};
\draw [fill] (13,4)  circle [radius=0.05];
\node at (13,3) [below] {};
\draw [fill] (13,3)  circle [radius=0.05]; 
\node at (13.5,-3) {}; 
\draw (10,0)--(11,1);
\draw (12,0)--(13,1);
\draw (11,1)--(12,2);
\draw (14,0)--(12,2);
\draw (12,2)--(13,3);
\draw (16,0)--(13,3);
\draw (13,3)--(13,4);
\end{tikzpicture}};
\node at (4, 3){\begin{tikzpicture}[scale=0.2] \node at (10,0) [below] {\scriptsize 1};
\draw [fill] (10,0)  circle [radius=0.05]; 
\node at (12,0) [below] {\scriptsize 4};
\draw [fill] (12,0)  circle [radius=0.05]; 
\node at (14,0) [below] {\scriptsize 2};
\draw [fill] (14,0)  circle [radius=0.05]; 
\node at (16,0) [below] {\scriptsize 3};
\draw [fill] (16,0)  circle [radius=0.05]; 
\node at (11,1) [below] {};
\draw [fill] (11,1)  circle [radius=0.05]; 
\node at (13,4) [left] {\scriptsize 5};
\draw [fill] (13,4)  circle [radius=0.05];
\node at (13,3) [below] {};
\draw [fill] (13,3)  circle [radius=0.05]; 
\node at (13.5,-3) {\scriptsize $F_j$}; 
\draw (10,0)--(11,1);
\draw (12,0)--(11,1);
\draw (11,1)--(12,2);
\draw (14,0)--(13,3);
\draw (12,2)--(13,3);
\draw (16,0)--(13,3);
\draw (13,4)--(13,3);
\end{tikzpicture}};
\node at (6, 3){\begin{tikzpicture}[scale=0.2] \node at (10,0) [below] {\scriptsize 1};
\draw [fill] (10,0)  circle [radius=0.05]; 
\node at (12,0) [below] {\scriptsize 2};
\draw [fill] (12,0)  circle [radius=0.05]; 
\node at (14,0) [below] {\scriptsize 4};
\draw [fill] (14,0)  circle [radius=0.05]; 
\node at (16,0) [below] {\scriptsize 3};
\draw [fill] (16,0)  circle [radius=0.05]; 
\node at (11,1) [below] {};
\node at (13,4) [left] {\scriptsize 5};
\draw [fill] (13,4)  circle [radius=0.05];
\node at (13,3) [below] {};
\draw [fill] (13,3)  circle [radius=0.05]; 
\node at (13.5,-3) {\scriptsize $F_k$}; 
\draw (10,0)--(11,1);
\draw (12,0)--(12,2);
\draw (11,1)--(12,2);
\draw (14,0)--(12,2);
\draw (12,2)--(13,3);
\draw (16,0)--(13,3);
\draw (13,3)--(13,4);
\end{tikzpicture}};
\node at (5,0) {\begin{tikzpicture}[scale=0.2] \node at (10,0) [below] {\scriptsize 1};
\draw [fill] (10,0)  circle [radius=0.05]; 
\node at (12,0) [below] {\scriptsize 4};
\draw [fill] (12,0)  circle [radius=0.05]; 
\node at (14,0) [below] {\scriptsize 2};
\draw [fill] (14,0)  circle [radius=0.05]; 
\node at (16,0) [below] {\scriptsize 3};
\draw [fill] (16,0)  circle [radius=0.05]; 
\node at (11,1) [below] {};
\draw [fill] (11,1)  circle [radius=0.05]; 
\node at (13,4) [left] {\scriptsize 5};
\draw [fill] (13,4)  circle [radius=0.05];
\node at (13,3) [below] {};
\draw [fill] (13,3)  circle [radius=0.05]; 
\node at (13.5,-3) {\scriptsize $F_{ijk}$}; 
\draw (10,0)--(11,1);
\draw (12,0)--(11,1);
\draw (11,1)--(12,2);
\draw (12,2)--(13,3);
\draw (16,0)--(13,3);
\draw (13,3)--(13,4);
\end{tikzpicture}};
\draw (5,1)--(4, 2);
\draw (5,1)--(6, 2);
\draw (5,1)--(8,2.25);
\draw (5,1)--(10,2.25);
\draw (5,1)--(2,2.25);
\draw (5,1)--(0,2.25);
\draw (4, 4)--(0, 5);
\draw (4, 4)--(2.5, 5); 
\draw (4, 4)--(5,5);
\draw (6, 4)--(5,5);
\draw (6, 4)--(7.5, 5);
\draw (6, 4)--(10, 5);
\draw (2,4)--(2.5,5);
\draw (0,4)--(0,5);
\draw (8,4)--(7.5,5);
\draw (10,4)--(10,5);
\node [above] at (5, 9) {$\hat{1}$};
\draw (5, 7)--(5,9);
\draw (0, 7)--(5,9);
\draw (2.5, 7)--(5,9);
\draw (7.5, 7)--(5,9);
\draw (10,7)--(5,9);
\node at (10,3){
\begin{tikzpicture}[scale=0.2] \node at (10,0) [below] {\scriptsize 1};
\draw [fill] (10,0)  circle [radius=0.05]; 
\node at (12,0) [below] {\scriptsize 2};
\draw [fill] (12,0)  circle [radius=0.05]; 
\node at (14,0) [below] {\scriptsize 4};
\draw [fill] (14,0)  circle [radius=0.05]; 
\node at (16,0) [below] {\scriptsize 3};
\draw [fill] (16,0)  circle [radius=0.05]; 
\node at (11,1) [below] {};
\node at (13,4) [left] {\scriptsize 5};
\draw [fill] (13,4)  circle [radius=0.05];
\node at (13,3) [below] {};
\draw [fill] (13,3)  circle [radius=0.05]; 
\node at (13.5,-3) {}; 
\draw (10,0)--(11,1);
\draw (12,0)--(10,0);
\draw (11,1)--(12,2);
\draw (14,0)--(12,2);
\draw (12,2)--(13,3);
\draw (16,0)--(13,3);
\draw (13,3)--(13,4);
\end{tikzpicture}};
\node at (8, 3){
\begin{tikzpicture}[scale=0.2] \node at (10,0) [below] {\scriptsize 1};
\draw [fill] (10,0)  circle [radius=0.05]; 
\node at (12,0) [below] {\scriptsize 2};
\draw [fill] (12,0)  circle [radius=0.05]; 
\node at (14,0) [below] {\scriptsize 4};
\draw [fill] (14,0)  circle [radius=0.05]; 
\node at (16,0) [below] {\scriptsize 3};
\draw [fill] (16,0)  circle [radius=0.05]; 
\node at (11,1) [below] {};
\node at (13,4) [left] {\scriptsize 5};
\draw [fill] (13,4)  circle [radius=0.05];
\node at (13,3) [below] {};
\draw [fill] (13,3)  circle [radius=0.05]; 
\node at (13.5,-3) {}; 
\draw (10,0)--(11,1);
\draw (12,0)--(14,0);
\draw (11,1)--(12,2);
\draw (14,0)--(12,2);
\draw (12,2)--(13,3);
\draw (16,0)--(13,3);
\draw (13,3)--(13,4);
\end{tikzpicture}};
\node at (2, 3){\begin{tikzpicture}[scale=0.2] \node at (10,0) [below] {\scriptsize 1};
\draw [fill] (10,0)  circle [radius=0.05]; 
\node at (12,0) [below] {\scriptsize 4};
\draw [fill] (12,0)  circle [radius=0.05]; 
\node at (14,0) [below] {\scriptsize 2};
\draw [fill] (14,0)  circle [radius=0.05]; 
\node at (16,0) [below] {\scriptsize 3};
\draw [fill] (16,0)  circle [radius=0.05]; 
\node at (11,1) [below] {};
\draw [fill] (11,1)  circle [radius=0.05]; 
\node at (13,4) [left] {\scriptsize 5};
\draw [fill] (13,4)  circle [radius=0.05];
\node at (13,3) [below] {};
\draw [fill] (13,3)  circle [radius=0.05]; 
\draw (10,0)--(11,1);
\draw (12,0)--(11,1);
\draw (11,1)--(12,2);
\draw (14,0)--(16,0);
\draw (12,2)--(13,3);
\draw (13,3)--(14,2);
\draw (16,0)--(13,3);
\draw (13,3)--(13,4);
\end{tikzpicture}};
\node at (0,3) {\begin{tikzpicture}[scale=0.2] \node at (10,0) [below] {\scriptsize 1};
\draw [fill] (10,0)  circle [radius=0.05]; 
\node at (12,0) [below] {\scriptsize 4};
\draw [fill] (12,0)  circle [radius=0.05]; 
\node at (14,0) [below] {\scriptsize 3};
\draw [fill] (14,0)  circle [radius=0.05]; 
\node at (16,0) [below] {\scriptsize 2};
\draw [fill] (16,0)  circle [radius=0.05]; 
\node at (11,1) [below] {};
\draw [fill] (11,1)  circle [radius=0.05]; 
\node at (13,3) [above] {\scriptsize 5};
\draw [fill] (13,3)  circle [radius=0.05];  
\draw (10,0)--(11,1);
\draw (12,0)--(11,1);
\draw (11,1)--(12,2);
\draw (14,0)--(12,2);
\draw (12,2)--(13,3);
\draw (16,0)--(13,3);
\end{tikzpicture}};
\end{tikzpicture}
\end{center}
\caption{An interval $[F_{ijk}, \hat{1}]$ in $\hat{S}([5])$ obtained as described in Definition \ref{problemdef}. The element $F_{ijk}$ is a maximal element of $I(C_j) \cap I(C_k)$, where $\{C_i, C_j, C_k\}$ is a critical triplet.}
\label{problempic}
\end{figure}

See Figure \ref{problempic} for an example of a critical triplet $\{C_i, C_j, C_k\}$ and the element $F_{ijk}$. 

\begin{figure}
\begin{center}
\begin{tikzpicture}
\node at (0,3){\begin{tikzpicture}[scale=0.3]
 \node at (10,0) [below] {\scriptsize 1};
\draw [fill] (10,0)  circle [radius=0.1]; 
\node at (12,0) [below] {\scriptsize 4};
\draw [fill] (12,0)  circle [radius=0.1]; 
\node at (14,0) [below] {\scriptsize 2};
\draw [fill] (14,0)  circle [radius=0.1]; 
\node at (16,0) [below] {\scriptsize 3};
\draw [fill] (16,0)  circle [radius=0.1]; 
\node at (11,1) [below] {};
\draw [fill] (11,1)  circle [radius=0.1]; 
\node at (13,4) [left] {\scriptsize 5};
\draw [fill] (13,4)  circle [radius=0.1];
\node at (13,3) [below] {};
\draw [fill, ] (13,3)  circle [radius=0.1]; 
\node at (13, -2) {\scriptsize A}; 
\draw (10,0)--(11,1);
\draw (12,0)--(11,1);
\draw (11,1)--(12,2);
\draw (12,2)--(13,3);
\draw (16,0)--(13,3);
\draw (13,3)--(13,4);
\draw[dashed] (14,0)--(13,3);
\end{tikzpicture}};
\node at (0,0) {\begin{tikzpicture}[scale=0.3]
 \node at (10,0) [below] {\scriptsize 1};
\draw [fill] (10,0)  circle [radius=0.1]; 
\node at (12,0) [below] {\scriptsize 4};
\draw [fill] (12,0)  circle [radius=0.1]; 
\node at (14,0) [below] {\scriptsize 2};
\draw [fill] (14,0)  circle [radius=0.1]; 
\node at (16,0) [below] {\scriptsize 3};
\draw [fill] (16,0)  circle [radius=0.1]; 
\node at (11,1) [below] {};
\draw [fill] (11,1)  circle [radius=0.1]; 
\node at (13,4) [left] {\scriptsize 5};
\draw [fill] (13,4)  circle [radius=0.1];
\node at (13,3) [below] {};
\draw [fill] (13,3)  circle [radius=0.1]; 
\node at (13,-2) {\scriptsize $F_{ijk}$}; 
\draw (10,0)--(11,1);
\draw (12,0)--(11,1);
\draw (11,1)--(12,2);
\draw (12,2)--(13,3);
\draw (16,0)--(13,3);
\draw (13,3)--(13,4);
\node [right] at (15,3) {\scriptsize $v(A)$};
\draw [->] (15,3)--(13.5,3);
\node [left] at (10.5,2.25) {\scriptsize $e(C)$};
\draw [->] (10.5,2.25) [out=60, in=135] to (12,2.25);
\end{tikzpicture}};
\draw (0,1)--(0,1.75);
\draw (0, 4)--(0,4.75);

\node at (0,6) {\begin{tikzpicture}[scale=0.3] \node at (10,0) [below] {\scriptsize 1};
\draw [fill] (10,0)  circle [radius=0.1]; 
\node at (12,0) [below] {\scriptsize 4};
\draw [fill] (12,0)  circle [radius=0.1]; 
\node at (14,0) [below] {\scriptsize 2};
\draw [fill] (14,0)  circle [radius=0.1]; 
\node at (16,0) [below] {\scriptsize 3};
\draw [fill] (16,0)  circle [radius=0.1]; 
\node at (11,1) [below] {};
\draw [fill] (11,1)  circle [radius=0.1]; 
\node at (13,4) [left] {\scriptsize 5};
\draw [fill] (13,4)  circle [radius=0.1];
\node at (13,3) [below] {};
\draw [fill] (13,3)  circle [radius=0.1]; 
\draw [fill] (12,2)  circle [radius=0.1]; 
\node at (13, -2) {\scriptsize C}; 
\draw (10,0)--(11,1);
\draw (12,0)--(11,1);
\draw (11,1)--(12,2);
\draw [dashed] (14,0)--(12,2);
\draw (12,2)--(13,3);
\draw (16,0)--(13,3);
\draw (13,3)--(13,4);
\end{tikzpicture}};
\end{tikzpicture}
\end{center}
\caption{Given the critical triplet $\{C_i, C_j, C_k\}$, every element covering $F_{ijk}$ can be obtained from $F_{ijk}$ by adding an edge between the singleton vertex and a vertex of $T(F_{ijk})$. Every element of $\mathcal{C}_{[n]}$ that is above $F_{ijk}$ can be obtained by adding a vertex to an edge of $T(F_{ijk})$, then adding an edge between this new vertex and the singleton vertex.}
\label{tfexample}
\end{figure}
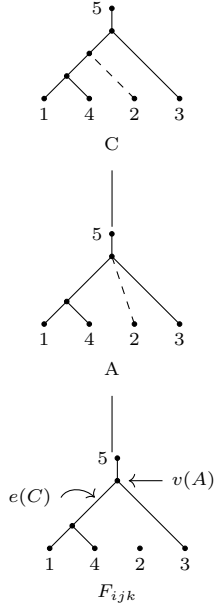

Given a critical triplet $\{C_i, C_j, C_k\}$, the element $F_{ijk}$ consists of an isolated vertex labeled $x$ and an $([n]-\{x\})$-tree. We will denote the $([n]-\{x\})$-tree of $F_{ijk}$ by $T(F_{ijk})$. If $A$ is an element of $\hat{S}([n])$ that covers $F_{ijk}$, then $A$ can be obtained from $F_{ijk}$ by adding an edge between the isolated vertex labeled $x$ and a unique vertex of $T(F_{ijk})$, say $v(A)$. We will thus associate the element $A$ to the vertex $v(A)$ of $T(F_{ijk})$. Similarly, if $C \in \mathcal{C}_{[n]}$ and $C$ is above $F_{ijk}$ in $\hat{S}([n])$, then $C$ can be obtained from $F_{ijk}$ by adding a new vertex $v$ to a unique edge $e(C)$ of $T(F_{ijk})$ and then adding an edge between $v$ and the isolated vertex labeled $x$. We will associate the coatom $C$ to the edge $e(C)$ of $T(F_{ijk})$. See Figure \ref{tfexample} for an example of a saturated chain $F_{ijk} \lessdot A \lessdot C$ of $\hat{S}([n])$ and to see how to obtain $A$ and $C$ from $T(F_{ijk})$.  In this way, the tree $T(F_{ijk})$ encodes information about all elements of $\hat{S}([n])$ that are above $F_{ijk}$. 

The next lemma shows that given an element $A$ covering $F_{ijk}$ in $\hat{S}([n])$, an element $C$ covers $A$ precisely when $e(C)$ is adjacent to $v(A)$ in $T(F_{ijk})$. 
\begin{lem} \label{tfstructure} 

Let $\{C_i, C_j, C_k\}$ be a critical triplet and let $C \in \mathcal{C}_{[n]}$ be such that $C > F_{ijk}$. Let $A \gtrdot F_{ijk} \in \hat{S}([n])$. 
 The element $C$ covers $A$ in $\hat{S}([n])$ if and only if $e(C)$ is adjacent to $v(A)$ in $T(F_{ijk})$.

\end{lem}
\begin{proof}
As noted above, $C$ can be obtained from $F_{ijk}$ by adding the vertex $v$ to $e(C)$ in $T(F_{ijk})$, then adding an edge between $v$ and the isolated vertex of $F_{ijk}$. Observe that given vertex $v(A)$ in $F_{ijk}$, $e(C)$ and $v(A)$ are adjacent in $T(F_{ijk})$ if and only if $v(A)$ and $v$ are adjacent in $C$. Contracting the edge between $v$ and $v(A)$ in $C$ yields $A$ if and only if $C$ covers $A$ in $\hat{S}([n])$. 
\end{proof}

Note that an immediate consequence of Lemma \ref{tfstructure} is that two elements of $\mathcal{C}_{[n]}$ cover a common element above $F_{ijk}$ in $\hat{S}([n])$ if and only if they correspond to edges that are adjacent to a common vertex in $T(F_{ijk})$.

Any $C \in \mathcal{C}_{[n]}$ that is above $F_{ijk}$ corresponds to both a vertex of $G_{NNI}(n)$ and an edge of $T(F_{ijk})$. We will need the following definition to help us connect $T(F_{ijk})$ and certain subgraphs of $G_{NNI}(n)$ in Corollary \ref{tfstructurecor}. 

\begin{defn} 
Let $F \in \hat{S}([n])$. An \textit{$F$-path} is a path in $G_{NNI}(n)$ all of whose vertices are elements of $\mathcal{C}_{[n]}$ that are above $F$. An \textit{$F$-cycle} is an $F$-path that is a simple cycle.
\end{defn}

If $C_i$, $C_j$, and $C_k$ are vertices that form a cycle in $G_{NNI}(n)$, Proposition \ref{3coatomscoverw} implies they are related by a NNI over some edge $\alpha$ and they cover a common element $w$. We will call such a cycle a \textit{trivial} cycle.  We will be most concerned with cycles of length greater than three in $G_{NNI}(n)$, which we will call \textit{nontrivial} cycles.

Lemma \ref{tfstructure} implies that $F_{ijk}$-paths in $G_{NNI}(n)$ correspond to sequences of adjacent edges in $T(F_{ijk})$. Since $T(F_{ijk})$ is a tree, we obtain the following corollary to this lemma:

\begin{cor} \label{tfstructurecor} Let $\{C_i, C_j, C_k\}$ be a critical triplet. 
\begin{enumerate}
\item There are no nontrivial $F_{ijk}$-cycles in $G_{NNI}(n)$. 
\item If for $C, C' \in \mathcal{C}_{[n]}$, $C$ and $C'$ are both above $F_{ijk}$, then there is an $F_{ijk}$-path between $C$ and $C'$ in $G_{NNI}(n)$. 
\end{enumerate}
\label{nonontrivcycles}
\label{ftree}
\end{cor}

The following proposition, which relies both on the fact that $T(F_{ijk})$ is a tree and that $F_{ijk}$-paths in $G_{NNI}(n)$ correspond to sequences of adjacent edges in $T(F_{ijk})$, is the last piece we need in order to prove Theorem \ref{s5norco}.
\begin{prop} If $\{C_i, C_j, C_k\}$ is a critical triplet, $C \in \mathcal{C}_{[n]}$ and $C$ is above $F_{ijk}$ in $\hat{S}([n])$, then an $F_{ijk}$-path from $C$ to $C_i$ can contain only one of $C_j$ or $C_k$. 
\label{fpathcjck}
\end{prop}
\begin{proof} Since $C_i$ covers both $F_j$ and $F_k$, $e(C_i)$ is the edge between $v(F_j)$ and $v(F_k)$ in $T(F_{ijk})$. By Lemma \ref{tfstructure}, $e(C_j)$ is also adjacent to $v(F_j)$ and $e(C_k)$ is also adjacent to $v(F_k)$. An $F_{ijk}$-path from $C$ to $C_i$ in $G_{NNI}(n)$ corresponds to a sequence of adjacent edges starting at $e(C)$ and ending at $e(C_i)$ in $T(F_{ijk})$. Without loss of generality, suppose  the $F_{ijk}$-path in $G_{NNI}(n)$ contains $C_j$, so it corresponds to a sequence of adjacent edges in $T(F_{ijk})$ that contains $e(C_j)$. Since $T(F_{ijk})$ is a tree, this sequence would stop at the edge $e(C_i)$ and thus never reach $e(C_k)$.

\end{proof}

We are now ready to prove our main theorem.

\begin{proof}[Proof of Theorem \ref{s5norco}]
We will prove that for $|X| \geq 5$, there does not exist a recursive coatom ordering of $\hat{S}(X)$; our result will then follow from Proposition \ref{shelliff}. 

Suppose that $\Omega$ is a recursive coatom ordering of $\hat{S}([n])$. Let $K$ be the nontrivial cycle in $G_{NNI}(n)$ shown in Figure \ref{3cycleex} and suppose $C_i$ is the vertex in $K$ that is latest in $\Omega$. Let $C_j$ and $C_k$ be the two vertices adjacent to $C_i$ in $K$. Observe that $\{C_i, C_j, C_k\}$ is a critical triplet, thus we can define the element $F_{ijk}$ as in Definition \ref{problemdef}. Let $A$ be the earliest coatom in $\Omega$ that is above $F_{ijk}$. 

The coatom $C_i$ is associated to the edge $e(C_i)$ in $T(F_{ijk})$; this is an internal edge since it is adjacent to $e(C_j)$ and $e(C_k)$ in $T(F_{ijk})$ and $e(C_j)$ and $e(C_k)$ are not adjacent by Proposition \ref{nonadjacent}. Since $e(C_i)$ is an internal edge of $T(F_{ijk})$, it is adjacent to four other edges in $T(F_{ijk})$. Thus $C_i$ is adjacent in $G_{NNI}(n)$ to exactly four other coatoms of $\hat{S}([n])$ that are above $F_{ijk}$. Two of these coatoms are $C_j$ and $C_k$. The other two we will call $D_j$ and $D_k$ where $D_j$ (respectively $D_k$) covers a common element with $C_j$ (resp. $C_k$) in $\hat{S}([n])$ and is adjacent to $C_j$ (resp. $C_k$) in $G_{NNI}(n)$. Consequently, any $F_{ijk}$-path from $A$ to $C_i$ contains at least one of $C_j$, $D_j$, $C_k$, and $D_k$, but since $C_j$ and $D_j$ (respectively $C_k$ and $D_k$) are adjacent in $G_{NNI}(n)$, if there is an $F_{ijk}$-path from $A$ to $C_i$ containing $D_j$ (resp. $D_k$), there is an $F_{ijk}$-path from $A$ to $C_i$ containing $C_j$ (resp. $C_k$). Without loss of generality, suppose there exists an $F_{ijk}$-path from $A$ to $C_i$ containing $C_k$; by Proposition \ref{fpathcjck}, it does not contain $C_j$. Let $B$ be the earliest coatom after $A$ in $\Omega$ that is above $F_{ijk}$ such that there exists an $F_{ijk}$-path between $B$ and $C_i$ containing $C_j$. This $B$ is guaranteed to exist: $C_j$ is in an $F_{ijk}$-path from $C_j$ to $C_i$, but if there is a coatom $C^*$ that comes before $C_j$ and after $A$ in $\Omega$ such that there exists an $F_{ijk}$-path from $C^*$ to $C_i$ containing $C_j$, then we take this earlier coatom $C^*$ to be $B$. Otherwise, we take $B=C_j$. Note $B$ necessarily comes before $C_i$ in $\Omega$.

Since $A$ and $B$ are both above $F_{ijk}$ in $\hat{S}([n])$, with $B$ coming after $A$ in $\Omega$, (ii) of Definition \ref{rco} implies the following: there must exist some element $C <_{\Omega} B$ and some element $z \in \hat{S}([n])$ such that $C \gtrdot z, B \gtrdot z$ and $z > F_{ijk}$ in $\hat{S}([n])$. We now show no such $z$ exists. Observe there is no $F_{ijk}$-path from $C$ to $C_i$ containing $C_j$, otherwise $C$ (and not $B$) would have been the first coatom after $A$ with the property that there is an $F_{ijk}$-path from it to $C_i$ containing $C_j$. We claim that any $C <_{\Omega} B$ is not adjacent to $B$ in $G_{NNI}(n)$, implying any $C <_{\Omega} B$ cannot cover a common element with $B$. To see this, note that if $C$ were adjacent to $B$, this would imply there was a nontrivial $F_{ijk}$-cycle in $G_{NNI}(n)$ consisting of the union of the $F_{ijk}$-path from $C_i$ to $B$ (containing $C_j$), the edge between $B$ and $C$, and the $F_{ijk}$-path from $C$ to $C_i$, contradicting Corollary \ref{nonontrivcycles}. Thus there is no $z>F_{ijk}$ covered by both $B$ and a coatom before $B$ in $\Omega$. We have arrived at a contradiction to the assumption that $\Omega$ is a recursive coatom ordering of $\hat{S}([n])$. 

Since there does not exist a recursive coatom ordering of $\hat{S}([n])$, Proposition \ref{shelliff} implies that for $|X| \geq 5$, the edge-product space $\mathcal{E}(X)$ is not shellable. 
\end{proof}

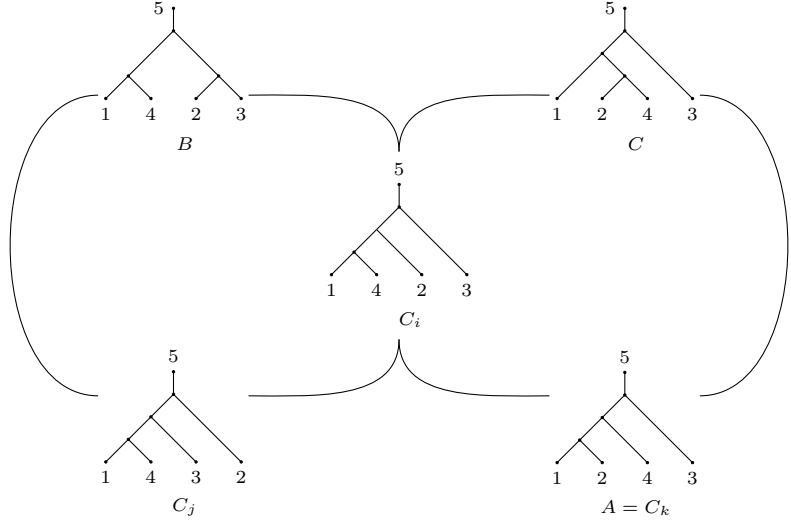
\begin{figure}
\begin{center}
\begin{tikzpicture}

\node at (1,-2.5){\begin{tikzpicture}[scale=0.3] 
\node at (10,0) [below] {\scriptsize 1};
\draw [fill] (10,0)  circle [radius=0.05]; 
\node at (12,0) [below] {\scriptsize 4};
\draw [fill] (12,0)  circle [radius=0.05]; 
\node at (14,0) [below] {\scriptsize 3};
\draw [fill] (14,0)  circle [radius=0.05]; 
\node at (16,0) [below] {\scriptsize 2};
\draw [fill] (16,0)  circle [radius=0.05]; 
\node at (11,1) [below] {};
\draw [fill] (11,1)  circle [radius=0.05]; 
\node at (13,4) [above] {\scriptsize 5};
\draw [fill] (13,4)  circle [radius=0.05];
\node at (13,3) [below] {};
\draw [fill] (13,3)  circle [radius=0.05]; 
\draw [fill] (12,2)  circle [radius=0.05]; 
\node at (13.5,-2) {\scriptsize $C_j$}; 
\draw (10,0)--(11,1);
\draw (12,0)--(11,1);
\draw (11,1)--(12,2);
\draw (14,0)--(12,2);
\draw (12,2)--(13,3);
\draw (16,0)--(13,3);
\draw (13,3)--(13,4);
\end{tikzpicture}};

\node at (4,0){\begin{tikzpicture}[scale=0.3] 
\node at (10,0) [below] {\scriptsize 1};
\draw [fill] (10,0)  circle [radius=0.05]; 
\node at (12,0) [below] {\scriptsize 4};
\draw [fill] (12,0)  circle [radius=0.05]; 
\node at (14,0) [below] {\scriptsize 2};
\draw [fill] (14,0)  circle [radius=0.05]; 
\node at (16,0) [below] {\scriptsize 3};
\draw [fill] (16,0)  circle [radius=0.05]; 
\node at (11,1) [below] {};
\draw [fill] (11,1)  circle [radius=0.05]; 
\node at (13,4) [above] {\scriptsize 5};
\draw [fill] (13,4)  circle [radius=0.05];
\node at (13,3) [below] {};
\draw [fill] (13,3)  circle [radius=0.05]; 
\node at (13.5,-2) {\scriptsize $C_i$}; 
\draw (10,0)--(11,1);
\draw (12,0)--(11,1);
\draw (11,1)--(12,2);
\draw (14,0)--(12,2);
\draw (12,2)--(13,3);
\draw (16,0)--(13,3);
\draw (13,3)--(13,4);
\end{tikzpicture}};

\node at (7,2.25){\begin{tikzpicture}[scale=0.3] 
\node at (10,0) [below] {\scriptsize 1};
\draw [fill] (10,0)  circle [radius=0.05]; 
\node at (12,0) [below] {\scriptsize 2};
\draw [fill] (12,0)  circle [radius=0.05]; 
\draw [fill] (13,1)  circle [radius=0.05]; 
\node at (14,0) [below] {\scriptsize 4};
\draw [fill] (14,0)  circle [radius=0.05]; 
\node at (16,0) [below] {\scriptsize 3};
\draw [fill] (16,0)  circle [radius=0.05]; 
\node at (11,1) [below] {}; 
\node at (13,4) [left] {\scriptsize 5};
\draw [fill] (13,4)  circle [radius=0.05];
\node at (13,3) [below] {};
\draw [fill] (13,3)  circle [radius=0.05]; 
\draw [fill] (12,2)  circle [radius=0.05]; 
\draw (10,0)--(11,1);
\draw (12,0)--(13,1);
\draw (11,1)--(12,2);
\draw (14,0)--(12,2);
\draw (12,2)--(13,3);
\draw (16,0)--(13,3);
\draw (13,3)--(13,4);
\node at (13.5,-2) {\scriptsize $C$}; 
\end{tikzpicture}};

\node at (7,-2.5){\begin{tikzpicture}[scale=0.3] 
\node at (10,0) [below] {\scriptsize 1};
\draw [fill] (10,0)  circle [radius=0.05]; 
\node at (12,0) [below] {\scriptsize 2};
\draw [fill] (12,0)  circle [radius=0.05]; 
\node at (14,0) [below] {\scriptsize 4};
\draw [fill] (14,0)  circle [radius=0.05]; 
\node at (16,0) [below] {\scriptsize 3};
\draw [fill] (16,0)  circle [radius=0.05]; 
\node at (11,1) [below] {};
\draw [fill] (11,1)  circle [radius=0.05]; 
\node at (13,4) [above] {\scriptsize 5};
\draw [fill] (13,4)  circle [radius=0.05];
\node at (13,3) [below] {};
\draw [fill] (12,2)  circle [radius=0.05]; 
\draw [fill] (13,3)  circle [radius=0.05]; 
\node at (13.5,-2) {\scriptsize $A=C_k$}; 
\draw (10,0)--(11,1);
\draw (12,0)--(11,1);
\draw (11,1)--(12,2);
\draw (14,0)--(12,2);
\draw (12,2)--(13,3);
\draw (16,0)--(13,3);
\draw (13,3)--(13,4);
\end{tikzpicture}};

\node at (1,2.25){\begin{tikzpicture}[scale=0.3] \node at (10,0) [below] {\scriptsize 1};
\draw [fill] (10,0)  circle [radius=0.05]; 
\node at (12,0) [below] {\scriptsize 4};
\draw [fill] (12,0)  circle [radius=0.05]; 
\node at (14,0) [below] {\scriptsize 2};
\draw [fill] (14,0)  circle [radius=0.05]; 
\node at (16,0) [below] {\scriptsize 3};
\draw [fill] (16,0)  circle [radius=0.05]; 
\node at (11,1) [below] {};
\draw [fill] (11,1)  circle [radius=0.05]; 
\node at (13,4) [left] {\scriptsize 5};
\draw [fill] (13,4)  circle [radius=0.05];
\node at (13,3) [below] {};
\draw [fill] (13,3)  circle [radius=0.05]; 
\draw [fill] (15,1)  circle [radius=0.05]; 
\draw (10,0)--(11,1);
\draw (12,0)--(11,1);
\draw (11,1)--(12,2);
\draw (14,0)--(15,1);
\draw (12,2)--(13,3);
\draw (16,0)--(13,3);
\draw (13,3)--(13,4);
\node at (13.5,-2) {\scriptsize $B$}; 
\end{tikzpicture}};

\draw (0,-2) to [out=180, in=180] (0,2);
\draw (2, 2) to [out=0, in=90] (4,1.25);
\draw (4,1.25) to [out=90, in=180] (6, 2);
\draw (8, 2) to [out=0, in=0] (8 ,-2);
\draw (2, -2) to [out=0, in=270] (4,-1.25);
\draw (4,-1.25) to [out=270, in=180] (6, -2);

\end{tikzpicture}
\end{center}
\caption{The subgraph of $G_{NNI}(5)$ consisting of coatoms of $\hat{S}([5])$ above $F_{ijk}$, obtained from the critical triplet $\{C_i, C_j, C_k\}$.}
\label{graphgf}
\end{figure}

\begin{ex} Here we provide an example of why a coatom ordering $\Omega$ of $\hat{S}([5])$ fails to be a recursive coatom ordering. Let $K$ be the cycle in $G_{NNI}(5)$ shown in Figure \ref{H5} with bold edges. Note that this is the same $K$ as in Figure \ref{3cycleex}.  Suppose that the element labeled $C_i$ in Figure \ref{graphgf} comes last in $\Omega$ among all elements of $K$. The element $C_i$ is adjacent to two other elements in $K$; we will call these two elements $C_j$ and $C_k$. The critical triplet $\{C_i, C_j, C_k\}$ is shown with the element $F_{ijk}$ in Figure \ref{problempic}. 

The graph shown in Figure \ref{graphgf} is the subgraph of $G_{NNI}(5)$ consisting of all the coatoms of $\hat{S}([5])$ that are above $F_{ijk}$. Let $A$ be the earliest coatom in $\Omega$ that is above $F_{ijk}$. As an example, suppose $A=C_k$. There is an $F_{ijk}$-path (i.e. a path in Figure \ref{graphgf}) from $A=C_k$ to $C_i$ that contains exactly one of $C_j$ or $C_k$. In our example, it contains $C_k$. Let $B$ be the first coatom after $A=C_k$ in $\Omega$ that is above $F_{ijk}$ and such that there exists a path between $B$ and $C_i$ containing $C_j$. As an example, suppose $B$ is the coatom labeled as such in Figure \ref{graphgf}. 

Both $A$ and $B$ are above $F_{ijk}$ in $\hat{S}([5])$. If $\Omega$ is a recursive coatom ordering, there must exist some coatom earlier than $B$ in $\Omega$, say $C$, and some element $z \in \hat{S}([5])$ such that $C \gtrdot z, B \gtrdot z$ and $z > F_{ijk}$. Based on how we were required to choose $B$, $C$ must be in an $F_{ijk}$-path from $A$ to $C_i$. In particular, $C$ must be the coatom labeled as such in \ref{graphgf}. However, $B$ and $C$ are not adjacent in $G_{NNI}(5)$, thus they cannot both cover an element $z$ in $\hat{S}([5])$. Thus no such $z$ and $C$ can exist. 
\label{s5ex}
\end{ex}

As noted earlier, the non-shellability of the edge-product space does not imply that the Tuffley poset is not shellable. However, the existence of the intervals $[F_{ijk}, \hat{1}]$ of $\hat{S}([n])$ described in the proof of Theorem \ref{s5norco} present an obstruction to the shellability of $S(X)$ as well. For similar reasons as those proving $\hat{S}([n])$ does not admit a recursive coatom ordering, in every ordering of the maximal chains of $S(X)$, there exists a chain $c$ so that the intersection of $c$ with the union of previous chains in the ordering is a chain of length two less than the length of $c$. This implies that in any ordering of facets of the order complex of $S(X)$, there exists a facet $F$ whose intersection with the union of earlier facets is not of dimension $dim(F)-1$. 

\bibliographystyle{model1-num-names}
\bibliography{nonshell}

\end{document}